\documentclass[a4paper,10pt]{amsart}
\usepackage{amsmath}
\usepackage[T1]{fontenc}
\usepackage{amssymb}
\usepackage{amsthm}
\newcommand{\md}{{\rm d}}

\newtheorem{thm}{Theorem}[section]
\newtheorem{prop}[thm]{Proposition}
\newtheorem{lem}[thm]{Lemma}
\newtheorem{corol}[thm]{Corollary}

\newtheorem{rem}[thm]{Remark}
\numberwithin{equation}{section}

\title[G-L functional with discontinuous constraint]
{Magnetic vortices for a Ginzburg-Landau type energy  with
discontinuous constraint}

\author{Ayman Kachmar}
\address{A. Kachmar\newline
Universit\'e Paris-Sud\\ D\'epartement de math\'ematique\\B\^at.
425\\F-91405 Orsay} \email{ayman.kachmar@math.u-psud.fr}
\subjclass[2000]{Primary 35J60; Secondary 35J20, 35J25, 35B40,
35Q55, 82D55} \keywords{Generalized Ginzburg-Landau energy
functional; proximity effects; global minimizers; unique positive
solution; vortices}

\date{November 22, 2007}

\begin{document}
\maketitle

\begin{abstract}
This paper is devoted to an analysis of vortex-nucleation
for a Ginzburg-Landau functional with
discontinuous constraint. This functional has been  proposed
as a model for vortex-pinning, and usually
accounts for  the energy
resulting from the interface of two superconductors. The
critical applied magnetic field for vortex nucleation is estimated in
the London singular limit,
and as a by-product,  results concerning vortex-pinning and
boundary conditions on the interface are obtained.
\end{abstract}

\tableofcontents
\section{Introduction and main results}
It is widely accepted among the physics community that spatial
inhomogeneities, impurities or point defects in a superconducting
sample provide pinning sites for vortices, preventing thus their
motion and the resultant induced resistivity, see \cite{Chetal, CR}
and the references therein. A similar behavior has also been
observed in superconducting samples subject to non-constant
temperatures, see \cite{MI}.\\
In the framework of the Ginzburg-Landau theory, it is proposed to
model the energy of an inhomogeneous superconducting sample by means
of the following functional (see \cite{Chetal, Ru})
\begin{equation}\label{V-EGL}
\mathcal G_{\varepsilon,H}(\psi,A)=\int_{\Omega}\left(
|(\nabla-iA)\psi|^2+\frac1{2\varepsilon^2}(p(x)-|\psi|^2)^2+ |{\rm
curl}\,A-H|^2\right)\,\md x.
\end{equation}
Here, $\Omega\subset\mathbb R^2$ is the 2-D cross section of the
superconducting sample, assumed to occupy a cylinder of infinite
height. The complex-valued function $\psi\in H^1(\Omega;\mathbb C)$
is called the `order parameter', whose modulus $|\psi|^2$ measures
the density of the superconducting electron Cooper pairs (hence
$\psi\equiv0$ corresponds to a normal state), and the real vector
field $ A=(A_1,A_2)$ is called the `magnetic potential', such that
the induced magnetic field in the
sample corresponds to ${\rm curl}\, A$.\\
The functional (\ref{V-EGL}) depends on many parameters:
$\frac1\varepsilon=\kappa$ is a characteristic of the
superconducting sample (a temperature independent quantity),
$H\geq0$ is the intensity of the applied magnetic field (assumed
constant and parallel to the axis of the superconducting sample),
$p(x)$ is a positive function modeling the impurities in the sample,
whose values are temperature dependent. The positive sign of the
function $p$ means that the temperature remains below the critical
temperature of the superconducting sample.\\
It is standard, starting from a minimizing sequence, to prove
existence of minimizers of the functional (\ref{V-EGL}) in the space
$H^1(\Omega;\mathbb C)\times H^1(\Omega;\mathbb R^2)$, see e.g.
\cite{GiPh}. A minimizer $(\psi,A)$ of (\ref{V-EGL}) is a weak
solution of the G-L equations:
\begin{equation}\label{G-L:equation}
\left\{
\begin{array}{l}
-(\nabla-iA)^2\psi=\displaystyle\frac1{\varepsilon^2}(p(x)-|\psi|^2)\psi,
\quad{\rm
 in}~\Omega,\\
\nabla^\bot\,{\rm
 curl}\,A=\big{(}i\psi,(\nabla-iA)\psi\big{)}\,\quad{\rm in}~\Omega,\\
n(x)\cdot (\nabla-iA)\psi=0,\quad{\rm curl}\,A=H\quad{\rm
  on}~\partial\Omega\,,
\end{array}\right.
\end{equation}
where $n(x)$ is the unite outward normal vector of
$\partial\Omega$.\\
It has been conjectured that for a minimizing configuration
$(\psi,A)$ of (\ref{V-EGL}), the vortices (zeros of $\psi$) should
be pinned near the minimal points of the function $p$ (or near the
critical points if $p$ is smooth), see \cite{CR, Ru}. Many authors
have addressed this question in the regime of extreme type II
superconducting materials, $\varepsilon\to0$. For instance,
Aftalion-Sandier-Serfaty \cite{AfSaSe} analyze the situation when
$p$ is periodic and smooth, Andr\'e-Baumann-Phillips \cite{AnBaPh}
analyze the situation when $p$ is smooth and having a finite number
of isolated zeros, and Alama-Bronsard \cite{AlBr} allow $p$ to have
negative values in some normal regions of the sample. We would also
like to mention the interesting work of Sigal-Ting \cite{ST}, who
prove existence and uniqueness of solutions with pinned vortices for
the Ginzburg-Landau
equation (\ref{G-L:equation}) in $\mathbb R^2$  when $H=0$ and the potential $p$ is in a suitable class.\\
In this paper, the function $p$ is a step function. We take
$\Omega=D(0,1)$ the unit disc in $\mathbb R^2$,  and
\begin{equation}\label{p(x)}
p(x)=\left\{\begin{array}{l} 1\quad{\rm if}~|x|\leq R\,,\\
a\quad{\rm if}~R<|x|<1\,,
\end{array}\right.
\end{equation}
where $a\in\mathbb R_+\setminus\{1\} $ and $0<R<1$ are given
constants.\\
Putting
\begin{equation}\label{V-S}
S_1=D(0,R),\quad S_2=D(0,1)\setminus \overline{D(0,R)}\,,
\end{equation}
then the above choice of $p$ has  two physical interpretations:
\begin{itemize}
\item $S_1$ and $S_2$ correspond to two  superconducting samples
with different critical temperatures;
\item The superconducting sample $\Omega$ is subject to two different
temperatures in the regions $S_1$ and $S_2$, which may happen by
cold or heat working $S_2$.
\end{itemize}
Lassoued-Mironescu analyze the functional (\ref{V-EGL}) without
magnetic field (i.e. $A=0$ \& $H=0$) and with $p$ as given in
(\ref{p(x)}), by assuming that minimizers satisfy a Dirichlet
boundary condition, $\psi=g$ on $\partial\Omega$, with $g$ valued in
$\mathbb S^1$ and has degree $d>0$, much in the same spirit of
B\'ethuel-Brezis-H\'elein \cite{BBH}. When $a>1$ and
$\varepsilon\to0$, they obtain that minimizers have $d$ vortices,
strictly localized in $S_1$, and whose positions are determined by a
finite dimensional problem (a
renormalized energy).\\
In this paper,
minimization of  the
functional (\ref{V-EGL}) will take place in  the space
$$\mathcal H=H^1(\Omega;\mathbb C)\times
H^1(\Omega;\mathbb R^2).$$ Thus we do not assume {\it a priori}
boundary conditions for admissible configurations, but minimizers
satisfy {\it natural boundary conditions}. We study nucleation of
vortices as the applied magnetic field varies, and we obtain that
their behavior is strongly dependent on the parameter $a$, leading
in some
situations (small values of $a$) to a {\it pinning} phenomenon.\\

We summarize in the next theorem some of the main results we have
obtained concerning the case of small values of $a$.

\begin{thm}\label{thm1}
Let   $(\psi_{\varepsilon,H},A_{\varepsilon,H})$
be a minimizer of (\ref{V-EGL}). There exists
a constant  $a_0\in]0,1[$, and for each $a\in]0,a_0[$, there
exist positive constants  $\lambda_*$, $\lambda_{\#}$, $\varepsilon_0$
and a function $]0,\varepsilon_0[\ni\varepsilon\mapsto
    k_\varepsilon\in\mathbb R_+$,
$0<\displaystyle\liminf_{\varepsilon\to0}k_\varepsilon\leq
\limsup_{\varepsilon\to0}k_\varepsilon<\infty$,
such that:
\begin{enumerate}
\item If $H<k_\varepsilon|\ln\varepsilon|
-\lambda_*\ln|\ln\varepsilon|$, then $|\psi_{\varepsilon,H}|\geq\displaystyle\frac{\sqrt{a}}2$ in
$\overline\Omega$.
\item If
  $H=k_\varepsilon|\ln\varepsilon|+\lambda\ln|\ln\varepsilon|$ and
  $\mu\geq-\mu_*$, then there exists a finite family of balls
$(B(a_i(\varepsilon),r_i(\varepsilon)))_{i\in I}$ with the following
properties:
\begin{enumerate}
\item
$\displaystyle\sum_{i\in I} r_i(\varepsilon)
<|\ln\varepsilon|^{-10}$\,;
\item $|\psi_{\varepsilon,H}|\geq\displaystyle\frac{\sqrt{a}}2$ in
$\Omega\setminus \displaystyle\bigcup_{i\in I}
  B(a_i(\varepsilon),r_i(\varepsilon))$\,;
\item
 Letting $d_i$ be the degree of
  $\psi_{\varepsilon,H}/|\psi_{\varepsilon,H}|$
on $\partial  B(a_i(\varepsilon),r_i(\varepsilon))$ if
  $B(a_i,r_i)\subset\Omega$ and $0$ otherwise, then we have
$$\displaystyle\sup_{\substack{i\in I\\|d_i|>0}} |R- |a_i(\varepsilon)|\,|\to
  0\quad{\rm  as}~
  \varepsilon\to0\,.$$
\item If $\lambda>\lambda_\#$
there exist
  positive constants $c$ and $C$ independent of $\varepsilon$ such that
$$c\ln|\ln\varepsilon|\leq \sum_i|d_i|\leq
  C\ln|\ln\varepsilon|\,\quad\forall~\varepsilon\in]0,\varepsilon_0]\,.$$
\end{enumerate}
\end{enumerate}
\end{thm}

Theorem~\ref{thm1} exhibits a completely different regime for the
nucleation of vortices when compared with the usual  G-L functional
defined in a simply connected domain \cite{SaSe}, and the result is
qualitatively  much more in the direction of a circular annulus
superconductor (c.f. \cite{AlBr}). In particular, Theorem~\ref{thm1}
states that vortices are localized near the less-superconducting
regions of the sample (i.e. $S_2$): This is the {\it pinning
phenomenon} predicted in the
physics literature, see e.g. \cite{CR}.\\

Let us define, as in \cite{SaSe}, the {\it vorticity measure}
\begin{equation}\label{vorticity}
\mu\left(\psi_{\varepsilon,H},A_{\varepsilon,H}\right)={\rm
curl}\left(i\psi_{\varepsilon,H},(\nabla-iA_{\varepsilon,H})\psi_{\varepsilon,H}\right)+{\rm
curl}\,A_{\varepsilon,H}\,.
\end{equation}
Theorem~\ref{thm1} shows that in the regime (d), we have (up to the
extraction of a subsequence)
$$\frac{\mu\left(\psi_{\varepsilon,H},A_{\varepsilon,H}\right)}{\ln|\ln\varepsilon|}\to\mu_*\quad{\rm
as}~\varepsilon\to0\,,$$ where $\mu_*$ is a measure supported in the
circle $\partial D(0,R)$. We conjecture that $\mu_*$ is indeed a
constant times the Lebesgue measure.\\

When $a>1$, we obtain a completely different behavior, which is that
of \cite{SaSe}. As in Theorem~\ref{thm1}, we get $k_\varepsilon>0$,
and a sequence of `critical fields'
$$H_{c,n}(\varepsilon)
=k_\varepsilon\left(|\ln\varepsilon|+(n-1)\ln|\ln\varepsilon|\right),\quad
(n\geq 1)\,,$$
such that, if $H_n<H<H_{n+1}$, then for a minimizer
$(\psi_{\varepsilon},
A_\varepsilon)$ of
(\ref{V-EGL}), $\psi_\varepsilon$ has exactly $n$ vortices
$\{x_i(\varepsilon)\}_{i=1}^n$, each of degree $1$, and there
positions are determined by minimizing a finite dimensional problem,
i.e.  {\it a renormalized energy} (see Section~\ref{Sec:a>1}).\\

\paragraph{\it Boundary conditions}\ \\
In addition to the pinning phenomenon, we obtain as a by-product
some interpretation concerning the boundary condition on the
$S_1$-$S_2$ interface. The precise result is the following.

\begin{thm}\label{thm2}
There exists a function $\mathbb
  R_+\setminus\{1\}\ni a\mapsto \gamma(a)\in\mathbb R\setminus\{0\}$
such that, if
$(\psi_{\varepsilon,H},A_{\varepsilon,H})$ is a minimizer of (\ref{V-EGL})
  satisfying $|\psi_{\varepsilon,H}|>0$ in $\overline\Omega$\,,
then the following limit holds:
\begin{equation}\label{deGennes}
\lim_{\varepsilon\to0}\,\varepsilon\left\|
n(x)\cdot(\nabla-iA_{\varepsilon,H})\psi_{\varepsilon,H}
+\frac{\gamma(a)}{\varepsilon}\psi_{\varepsilon,H}\right\|_{L^2(\partial
D(0,R))}=0\,,
\end{equation}
where $n(x)=\displaystyle\frac{x}{|x|}$
for all $x\in\mathbb
R^2\setminus\{0\}$, is the outward unit normal vector.\\
Furthermore, the function $\gamma$ satisfies: (1) $\gamma(a)>0$ if
$a<1$; (2) $\gamma(a)<0$ if $a>1$.
\end{thm}

Thus, below the first critical field $H_{C_1}$, minimizers
approximately  satisfy a Robin-type boundary condition on the $S_1$-boundary:
\begin{equation}\label{Robin}
n(x)\cdot (\nabla-iA_{\varepsilon,H})\psi_{\varepsilon,H}
+\frac{\gamma(a)}\varepsilon\psi_{\varepsilon,H}(1+o(1))=0 \quad {\rm
  on}~\partial S_1.\end{equation}
This is a boundary condition of the type predicted by de\,Gennes
($\gamma(a)$
being called {\it the de\,Gennes parameter}), see
\cite{deGe}. When $a>1$, $\gamma(a)<0$, hence we justify the modeling
of Fink-Joiner \cite{FJ}, who use
a negative `de\,Gennes parameter' to model a superconductor
surrounded by another superconductor with a higher critical
temperature. They claim also that this is the setting when cold
working the surface of
superconducting samples (see \cite{ICM, MI} for more recent reviews of
this topic). \\
The result of Theorem~\ref{thm2} also justifies the analysis we carried
out in \cite{kach5, kach2, kach1} for problems involving boundary conditions
of the  type (\ref{Robin}),
and complements results in this direction obtained in \cite{kach4, kach3}.\\
\paragraph{\it Main points of the proofs.}\ \\
Let us briefly describe the main points of the proofs
of Theorems~\ref{thm1}
and \ref{thm2}, and thus explain what stands behind their
statements.\\
The starting point is the analysis of minimizers of (\ref{V-EGL}) when
$H=0$. In this case, (\ref{V-EGL}) has, up to a gauge transformation,
a unique minimizer $(u_\varepsilon,0)$ where $u_\varepsilon$ is a
positive real-valued function. The asymptotic profile  of
$u_\varepsilon$
as $\varepsilon\to0$ is obtained in Theorem~\ref{mainthm-H=0}, which
proves Theorem~\ref{thm2} when $H=0$ with
a stronger convergence in $L^\infty$ norm.\\
When $H>0$, let $(\psi,A)$ be a minimizer of (\ref{V-EGL}).
Inspired by Lassoued-Mironescu \cite{LaMi},
we introduce a {\it normalized density}\footnote{Notice that $\varphi$
  and $\psi$ have the same vortices.}
$$\varphi=\frac{\psi}{u_\varepsilon}\,.$$
Then $|\varphi|\leq1$ and we are led to the analysis of the following
functional
(see Lemma~\ref{V-lem-psi<u})
\begin{equation}\label{V-EGL*}
\mathcal F_{\varepsilon,H}(\varphi,A)=\int_\Omega
\left(u_\varepsilon^2
|(\nabla-iA)\varphi|^2
+\frac1{2\varepsilon^2}
u_\varepsilon^4(1-|\varphi|^2)^2+ |{\rm
curl}\,A-H|^2\right)\md x\,,\end{equation}
using tools from Sandier-Serfaty \cite{SaSe}.\\
When we take $\varphi=1$ in
(\ref{V-EGL*}) and we minimize the resulting functional over $A\in
H^1(\Omega;\mathbb R^2)$, we get that the minimizer is
$\frac{H}{u_\varepsilon^2}\nabla^\bot h_\varepsilon$, where
$h_\varepsilon:\Omega\longrightarrow]0,1[$ is the solution of the
equation:
\begin{equation}\label{London*}
-{\rm div}\left(\frac1{u_\varepsilon^2}\nabla
h_\varepsilon\right)+h_\varepsilon=0\quad{\rm in}~\Omega\,,\quad
h_\varepsilon=1\quad{\rm on}~\partial\Omega\,.
\end{equation}
The constant $k_\varepsilon$ appearing in Theorem~\ref{thm1} is defined by
\begin{equation}\label{k-eps*}
k_\varepsilon=
\displaystyle
\frac12\left(\max_{x\in\overline\Omega}
\frac{1-h_\varepsilon(x)}{u_\varepsilon^2(x)}\right)^{-1}.\end{equation}
Thanks to our choice of the domain $\Omega$ and the step function
$p(x)$ in (\ref{p(x)}),
we show that the function $h_\varepsilon(x)$ is radially
symmetric and strictly increasing with respect to $|x|$
(see Lemma~\ref{V-lem-hepsilon}). This
permits to show that
$$0<\liminf_{\varepsilon\to0}k_\varepsilon\leq\limsup_{\varepsilon\to0}k_\varepsilon<+\infty\,.$$
Roughly speaking, the analysis of Sandier-Serfaty (c.f. \cite{SaSe})
says that near the first critical magnetic field, the vortices of a
minimizer of\footnote{These are also the vortices of a minimizer
of (\ref{V-EGL}).} (\ref{V-EGL*})  are localized as $\varepsilon\to0$ near
the set $\Lambda_\varepsilon=\left\{x\in\overline\Omega~:~
\displaystyle\frac{1-h_\varepsilon(x)}{u_\varepsilon^2(x)}=\frac12k_\varepsilon^{-1}\right\}$.
We localize the set $\Lambda_\varepsilon$ by means of a fine
semi-classical analysis. We obtain when $a$ is sufficiently small that
the set $\Lambda_\varepsilon$ consists of a circle $\partial
B(0,R_\varepsilon)$, where $R_\varepsilon\in]R,1[$ has the following
asymptotic behavior (see Theorem~\ref{corollary}):
$$\varepsilon\ll R_\varepsilon-R\ll \varepsilon^\alpha\quad{\rm as~}\varepsilon\to0\,,\quad
\quad\quad (\alpha\in]0,1[ {\rm~is~given}).$$ Let us mention that when
$a>1$,  we prove that the set $\Lambda_\varepsilon$ consists of a
    single point, $\Lambda_\varepsilon=\{0\}$, and for this reason,
minimizers of (\ref{V-EGL}) exhibit the same behavior as the
one present in \cite{SaSe}, i.e. near the first critical magnetic field,
a minimizer has  a finite number of vortices localized near the
center of the disc and whose exact positions are determined by a
finite dimensional problem  (a renormalized energy).\\

\paragraph{\it Outline of the paper.}\ \\
Section~\ref{V-Sec-EnH=0} is devoted to a preliminary analysis of the
variational problem (\ref{V-EGL}). In particular, a detailed analysis
is given for the case without magnetic field $H=0$.\\
Section~\ref{V-Sec-meissnerstate} is devoted to an analysis of
the equation (\ref{London*}).\\
Section~\ref{Section:LB} is devoted to derive a lower bound of the
minimizing energy, involving the construction of vortex-balls.\\
Section~\ref{Section:UB} is devoted to establish an upper bound of the
minimizing energy, that is involved with a careful analysis of a
Green's potential.\\
Finally, Section~\ref{section:proofs} is devoted to the proofs of
Theorems~\ref{thm1} and \ref{thm2}, through the matching of the lower
and upper bounds obtained in Sections~\ref{Section:LB} and
\ref{Section:UB}
respectively.

\paragraph{\it A remark on the notation.}\ \\
The letters $C,\widetilde C, M,$ etc. will  denote positive constants
independent of $\varepsilon$. For $n\in\mathbb N$ and $X\subset\mathbb
R^n$,
$|X|$ denotes the Lebesgue measure of $X$. $B(x,r)$ denotes the open
ball in $\mathbb R^n$ of radius $r$ and center $x$. $(\cdot,\cdot)$ denotes the
scalar product in $\mathbb C$ when identified with $\mathbb R^2$.
For two positive functions
$a(\varepsilon)$ and $b(\varepsilon)$, we write
$a(\varepsilon)\ll b(\varepsilon)$ as $\varepsilon\to0$ to mean that
$\displaystyle\lim_{\varepsilon\to0}\frac{a(\varepsilon)}{b(\varepsilon)}=0$.

\section{Preliminary analysis of minimizers}
\label{V-Sec-EnH=0}

\subsection{The case without applied magnetic field}
This section is devoted to an analysis for minimizers of (\ref{V-EGL}) when the
applied magnetic field $H=0$. We follow closely
similar results obtained in
\cite{kach3}.\\
We keep the notation introduced in Section~1. Upon taking $A=0$ and
$H=0$ in (\ref{V-EGL}), one is led to introduce the functional
\begin{equation}\label{V-EnH=0}
\mathcal G_\varepsilon(u):=\int_\Omega\left(|\nabla u|^2
+\frac{1}{2\varepsilon^2}(p(x)-u^2)^2\right)\,\md
x\,,
\end{equation}
defined for functions in $H^1(\Omega;\mathbb R)$.\\
We introduce
\begin{equation}\label{V-C0}
C_0(\varepsilon)=\inf_{u\in H^1(\Omega;\mathbb R)}
\mathcal G_\varepsilon(u)\,.
\end{equation}

The next theorem is an analogue of Theorem~1.1 in \cite{kach3}.

\begin{thm}\label{V-thm-kach3}
Given $a\in\mathbb R_+\setminus\{1\}$, there exists $\varepsilon_0$ such that for all
$\varepsilon\in]0,\varepsilon_0[$, the functional (\ref{V-EnH=0})
admits in $H^1(\Omega;\mathbb R)$ a
minimizer $u_\varepsilon\in C^2(\overline{S_1})\cup C^2(\overline{S_2})$
such that
$$\min(1,\sqrt{a}\,)<u_\varepsilon<\max(1,\sqrt{a}\,)\quad {\rm in}~\overline{\Omega}.$$
Furthermore, with our choice of the domains $\Omega,S_1$ and $S_2$ in
(\ref{V-S}), the function $u_\varepsilon$ is radial.\\
If $H=0$, minimizers of (\ref{V-EGL}) are gauge equivalent to the
state $(u_\varepsilon,0)$.
\end{thm}

The asymptotic behaviour of the function $u_\varepsilon$ when
$\varepsilon\to0$ is based on the understanding of the following
canonical equation:
\begin{equation}\label{Canon-equation}
\left\{
\begin{array}{l}
-\Delta u=(1-u^2)u~{\rm in}~\mathbb R\times\mathbb R_-,\quad
-\Delta u=(a-u^2)u~{\rm in}~\mathbb R\times\mathbb R_+,\\
\partial _{x_2}u(\cdot,0_-)=\partial _{x_2}u(\cdot,0_+),\quad
u(\cdot,0_-)=u(\cdot,0_+)\quad{\rm on}~\mathbb R.
\end{array}\right.
\end{equation}
When $a\not=1$, it is easy to verify that (\ref{Canon-equation}) has
the following solution
\begin{equation}\label{canon-sol}
\mathbb R^2\ni(x_1,x_2)\mapsto U(x_2)\,,
\end{equation}
where the function $U(x_2)$ is defined by
\begin{equation}\label{U(x_2)}
U(x_2)=\left\{
\begin{array}{l}
\displaystyle\frac{\beta_1(a)\,\exp(-\sqrt{2}\,x_2)-1}
{\beta_1(a)\,\exp(-\sqrt{2}\,x_2)+1}\,,\quad{\rm if}~x_2\in\mathbb
R_-\,,\\
\sqrt{a}\,
\displaystyle\frac{\beta_2(a)\,\exp(\sqrt{2/a}\,x_2)-1}
{\beta_2(a)\,\exp(\sqrt{2/a}\,x_2)+1}\,,\quad{\rm if}~x_2\in\mathbb
R_+\,.\end{array}\right.\end{equation}
The constants $\beta_1(a)$ and $\beta_2(a)$ are given explicitly:
\begin{equation}\label{beta}
\beta_1(a)=\frac{\alpha(1+\alpha\sqrt{a})}{\alpha-\sqrt{a}},\quad
\beta_2(a)=-\alpha^2\beta_1(a),\quad
\alpha=\frac{1+\sqrt{a}-\sqrt{2(1+a)}}{1-\sqrt{a}}\,.
\end{equation}
Furthermore, we have the following properties:
\begin{equation}\label{properties1}
\left\{
\begin{array}{l}
\forall~a\in]0,1[\,,\quad\beta_1(a)>1~\&~\beta_2(a)<-1\,;\\
\forall~a\in]1,\infty[\,,\quad\beta_1(a)<-1~\&~\beta_2(a)>1,
\end{array}\right.
\end{equation}
and
\begin{equation}
U'(0)=\gamma(a)\,U(0),\quad\gamma(a)=\alpha
\frac{a\alpha^3+\sqrt{a}\,\alpha^2+a\alpha+\sqrt{a}}{\alpha^3+
(4-\sqrt{a})\alpha^2-3\sqrt{a}\,\alpha+a}\,.\label{properties2}
\end{equation}

As in \cite[Theorem~1.5]{kach3}, we get that the solution given by
(\ref{canon-sol}) is unique in a certain class of functions.

\begin{thm}\label{thm-canonEq}
Let $a\in\mathbb R_+\setminus\{1\}$. Eq. (\ref{Canon-equation})
admits, in the class of functions
$\mathcal C=\{u\in H^2_{\rm loc}(\mathbb R^2)\cap L^\infty(\mathbb
R^2)~:~u\geq 0~{\rm in}~\mathbb R^2\}$,
a unique non-trivial solution given by
(\ref{canon-sol}).\end{thm}
\paragraph{\bf Proof.}
Since the proof is very close to that of \cite[Theorem~1.5]{kach3},
we sketch only the main steps.\\
By adjusting the proof of \cite[Lemma~4.2]{kach3}, we obtain that if
$u\not\equiv0$ solves (\ref{Canon-equation}), then $0<u<1$ in $\mathbb
R^2$. This permits us, when following step by step the proof of
\cite[Lemma~4.3]{kach3} and \cite[Lemma~5.3]{LuPa96},
to get a positive constant $C\in]0,1[$
such that for any solution $u$ of (\ref{Canon-equation}) in $\mathcal C$, we
have
\begin{equation}\label{kach3-4.9}
\inf_{x\in\mathbb R^2}u(x)>C\,.\end{equation}
Also, we prove in \cite[Lemma~4.4]{kach3} that, for $u\in\mathcal C$
a solution of (\ref{Canon-equation}),
\begin{equation}\label{kach3-4.11}
\lim_{x_2\to-\infty}\left(\sup_{x_1\in\mathbb
  R}(1-u(x_1,x_2))\right)=0\,,\quad
\lim_{x_2\to+\infty}\left(\sup_{x_1\in\mathbb
  R}(\sqrt{a}-u(x_1,x_2))\right)=0\,.
\end{equation}
Now, let $u_1,u_2\in\mathcal C$ be solutions of (\ref{Canon-equation}). We
introduce
\begin{equation}\label{kach3-H-lambda}
\lambda_*=\sup\{\lambda\in[0,1[~:~u_2(x)> \lambda u_1(x)\}\,.
\end{equation}
Then, by (\ref{kach3-4.9}), $\lambda_*>0$. We claim that
$\lambda_*=1$. Once this is shown to hold,
Theorem~\ref{thm-canonEq} is proved.\\
We argue by contradiction: If $\lambda_*<1$, then
\begin{equation}\label{kach3-w}
\inf_{x\in\mathbb R^2} w(x)=0\,,
\end{equation}
where $w(x)=u_2(x)-\lambda_*u_1(x)$. Now, let
$(x_n)=\left((x_n^1,x_n^2)\right)$ be
a minimizing sequence:
$$\lim_{n\to+\infty}w(x_n)=0\,.$$
Since the maximum principle yields that $w(x)>0$ for all $x$, the
sequence $(x_n)$ should be unbounded, hence we assume that
$\lim_{n\to+\infty }|x_n|=+\infty$. Also, by (\ref{kach3-4.11}), $(x_n^2)$
should be bounded, hence we assume that $\lim_{n\to+\infty}x_n^2=b$.\\
Now, the functions $u_j^n(x_1,x_2)=u_j(x_1+x_1^n,x_2)$, $j=1,2$, solve
(\ref{Canon-equation}) in $\mathcal C$, and up to extraction of a
subsequence, they  converge locally to functions\break
$\widetilde u_j\in C_{\rm  loc}^2(\overline{\mathbb R\times\mathbb
R_\pm}\,;\mathbb R)$, $j=1,2$. Now, $\widetilde u_1$, $\widetilde u_2$
solve (\ref{Canon-equation}) in
$\mathcal C$, $\widetilde u_2\geq \lambda_* \widetilde u_1$ and
$\widetilde u_2(0,b)=\lambda_* \widetilde u_1(0,b)$. On the other hand, the
strong maximum principle insures that $\widetilde u_2(x)>\lambda_*
\widetilde u_1(x)$ for all $x\in\mathbb R^2$, hence we have a
contradiction.\hfill$\Box$

\begin{rem}\label{a=1}
It is known (see the remark p.~163 in \cite{LuPa96})
that when $a=1$, the trivial solution $U\equiv1$ is the unique
positive and bounded solution of Eq. (\ref{Canon-equation}).
\end{rem}

By a blow-up argument, Theorem~\ref{thm-canonEq} permits us to obtain
the asymptotic behaviour of the minimizer $u_\varepsilon$ of
(\ref{V-EnH=0}).

\begin{thm}\label{mainthm-H=0}
Let $a\in\mathbb R_+\setminus\{1\}$ and $u_\varepsilon$
be the positive minimizer of (\ref{V-EnH=0})
introduced in Theorem~\ref{V-thm-kach3}. Then, the following
asymptotics hold as $\varepsilon\to0$~:
\begin{equation}\label{asymptot1-H=0}
\lim_{\varepsilon\to0}\left\|u_\varepsilon(x)
-U\left(\frac{|x|-R}\varepsilon\right)\right\|_{L^\infty(\Omega)}=0\,,
\end{equation}
\begin{equation}\label{asymptot2-H=0}
\forall~C>0,\quad
\lim_{\varepsilon\to0}\varepsilon\left\|u_\varepsilon(x)
-U\left(\frac{|x|-R}\varepsilon\right)\right\|
_{W^{1,\infty}(\{x\in\mathbb R^2:|R-|x||\leq C\varepsilon\})}=0\,,
\end{equation}
where $U$ is the function introduced in (\ref{canon-sol}).
\end{thm}

In particular, Theorem~\ref{mainthm-H=0} provides a stronger version of
Theorem~\ref{thm2} when $H=0$.\\

\paragraph{\bf Proof of Theorem~\ref{mainthm-H=0}.}
Let $(r,\theta)$ be polar coordinates, $0<r<1$,
$-\pi\leq\theta<\pi$,  and set
$$t=r-R,\quad s=R\,\theta\,.$$
Given $s_0\in[-R\pi,R\pi[$, we define the rescaled function,
$$\widetilde u_\varepsilon(s,t)=u_\varepsilon\left((R+\varepsilon t)
e^{i\varepsilon
(s-s_0)/R}\right)\,,\quad  \frac{R-1}\varepsilon<t<\frac{1-R}\varepsilon,
\quad -\pi\frac{R}\varepsilon<s-s_0<\pi\frac{R}\varepsilon.$$
The equation of $\widetilde u_\varepsilon$ becomes:
$$\left\{
\begin{array}{l}
-\Delta_\varepsilon\, \widetilde u_\varepsilon =(1-\widetilde
u_\varepsilon^2)\widetilde u_\varepsilon,
\quad \frac{R-1}\varepsilon<t<0,~|s-s_0|<\pi\frac{R}\varepsilon,
\\ \\
-\Delta_\varepsilon\,\widetilde u_\varepsilon=(a-\widetilde
u_\varepsilon^2)
\widetilde
u_\varepsilon,
\quad 0<t<\frac{1-R}\varepsilon,~|s-s_0|<\pi\frac{R}\varepsilon,\\ \\
\displaystyle\frac{\partial\widetilde u_\varepsilon}{\partial
t}(\cdot,0_-)=\displaystyle\frac{\partial\widetilde
u_\varepsilon}{\partial t}(\cdot,0_+),\quad
\widetilde u_\varepsilon(\cdot,0_-)=\widetilde u(\cdot,0_+)\,,
\end{array}\right.$$
where
$$\Delta_\varepsilon=\left(1-\varepsilon
\frac{t}R\right)^{-2}\partial_s^2
+\partial_t^2 -\frac{\varepsilon }{\left(R-\varepsilon t\right)}\partial_t.$$
Now, by elliptic estimates, the function $\widetilde u_\varepsilon$
converges to a function $u$ in $W^{2,\infty}_{\rm loc}(\mathbb
R^2)$. Furthermore, $u$ solves (\ref{Canon-equation}) in $\mathcal C$, and
by \cite[Lemma~5.2]{kach3}, there exist constants $k_0,c_0>0$ such
that $u(0,k_0)>c_0$.
Thus, we conclude by Theorem~\ref{thm-canonEq} that $u(s,t)=
U(t)$, where $U$ is given in (\ref{canon-sol}), and therefore, by
coming back to the initial scale,
\begin{equation}\label{convergence}
\forall~C>0,\quad\forall~k\in\{0,1,2\},\quad\lim_{\varepsilon\to 0}
\varepsilon^k\left\|u_\varepsilon(s,t)
-U\left(\frac t\varepsilon\right)\right\|_{W^{k,\infty}(\{|s-s_0|\leq
C\varepsilon,\,|t|\leq C\varepsilon\})}=0.\end{equation}
To prove (\ref{asymptot1-H=0}), let
$x_\varepsilon=(R+t(x_\varepsilon))
e^{i\,s(x_\varepsilon)/R}\in\overline\Omega$ such that
$$\left|u_\varepsilon(x_\varepsilon)
-U\left(\frac{|x_\varepsilon|-R}\varepsilon\right)\right|
=\left\|u_\varepsilon(x)
-U\left(\frac{|x|-R}\varepsilon\right)\right\|_{L^\infty(\Omega)}.$$
If $|R-|x_\varepsilon|\,|/\varepsilon$ is bounded,
then (\ref{asymptot1-H=0}) becomes a consequence of
(\ref{convergence}) upon taking $s_0=s(x_\varepsilon)$. Otherwise, if
$\displaystyle\lim_{\varepsilon\to0}
\frac{R-|x_\varepsilon|}\varepsilon=\pm\infty$,
we get again by a blow-up argument that
$u_\varepsilon(x_\varepsilon)\to 1$ if the limit is $+\infty$, and
$u_\varepsilon(x_\varepsilon)\to \sqrt{a}$ if the limit is $-\infty$. This
establishes (\ref{asymptot1-H=0}) in this case.\\
The asymptotic limit (\ref{asymptot2-H=0}) is also a simple
consequence of (\ref{convergence}). We take
$y_\varepsilon=(R+t(y_\varepsilon))e^{i\,s(y_\varepsilon)/R}$ such that
$$
\left\|\nabla\left(u_\varepsilon(x)
-U\left(\frac{|x|-R}\varepsilon\right)\right)\right\|_
{L^\infty(|R-|x||\leq C\varepsilon\})}=
\left|\nabla\left(u_\varepsilon(y_\varepsilon)
-U\left(\frac{|y_\varepsilon|-R}\varepsilon\right)\right)\right|\,.$$
Then we apply (\ref{convergence}) with $s_0=s(y_\varepsilon)$.
\hfill$\Box$\\

We state also some estimates, taken from \cite[Proposition~5.1]{kach3},
that describe the decay of  $u_\varepsilon$ away from the
boundary of $S_1$.

\begin{lem}\label{interior}
Let $k\in\mathbb N$.
There exist positive constants $\varepsilon_0$, $\delta$ and $C$ such
that, for all $\varepsilon\in]0,\varepsilon_0]$,
\begin{equation}\label{interiordecay}
\left\|(1-u_\varepsilon)\exp\left(\frac{\delta|R-|x|\,|}\varepsilon\right)\right\|_{H^k(S_1)}+\left\|(\sqrt{a}-u_\varepsilon)\exp\left(\frac{\delta|R-|x|\,|}\varepsilon\right)\right\|_{H^k(S_2)}\leq\frac{C}{\varepsilon^k}\,.
\end{equation}
\end{lem}

Another property that we need is the monotonicity of the function
$u_\varepsilon$ (recall that, in the setting of
Theorem~\ref{mainthm-H=0}, $u_\varepsilon$ is radial).

\begin{lem}\label{monotonicity}
With the choice of $\Omega$, $S_1$ and $S_2$ as in (\ref{V-S}), the
function $u_\varepsilon$ is increasing if $a>1$ and decreasing
if $a<1$.
\end{lem}
\paragraph{\bf Proof.}
We only prove the result of the lemma  when $a<1$.\\
{\it Step~1. $u'_\varepsilon(R)\not=0$.}\\
Notice that $u_\varepsilon$ is positive and satisfies the equations:
\begin{eqnarray}
&&-u''_\varepsilon-\frac1ru'_\varepsilon=\frac1{\varepsilon^2}(1-u_\varepsilon^2)u_\varepsilon\quad{\rm
    in}~]0,R[\label{u-Eq1}\\
&&-u''_\varepsilon-\frac1ru'_\varepsilon=
\frac1{\varepsilon^2}(a-u_\varepsilon^2)u_\varepsilon\quad{\rm
    in}~]R,1[
\label{u-Eq2}\\
&&u'_\varepsilon(0)=0\,,\quad u'_\varepsilon(1)=0\,.\label{u-Eq3}
\end{eqnarray}
Therefore, if $u'_\varepsilon(R)=0$, then
$$u_\varepsilon\equiv 1\quad{\rm in}~S_1\,,\quad
u_\varepsilon\equiv a\quad {\rm in}~S_2\,.$$
This is impossible since the function $u_\varepsilon$ is
in $H^1(\Omega)$.\\
{\it Step~2. The function $u_\varepsilon$ is decreasing in $[0,R]$.}\\
Recall that, by Theorem~\ref{V-thm-kach3}, $\sqrt{a}\,<
u_\varepsilon< 1$ in $\overline\Omega$. It is then easy to verify
from Eqs. (\ref{u-Eq1}) and (\ref{u-Eq3}) that
$u''_\varepsilon(0)>0$. Let us denote by $\widetilde u_\varepsilon$
the even extension of $u_\varepsilon$ in $]-R,0[$. Then it is easy
to verify that (i) $\widetilde u_\varepsilon\in C^2([-R,R])$\,; (ii)
If $r_0\in]-R,R[$ is a critical point of $\widetilde
    u_\varepsilon$, then $\widetilde u_\varepsilon''(r_0)<0$\,.
This shows that every critical point of $\widetilde u_\varepsilon$
is a local maximum. Consequently, $\widetilde u_\varepsilon$ should
have a unique critical point in $]-R,R[$ and $\widetilde
u_\varepsilon'$ should change its sign only in this critical point.
Since $\widetilde u'_\varepsilon(0)=0$ and $\widetilde
u_\varepsilon''(0)<0$, we deduce that $\widetilde u'_\varepsilon< 0$
in  $]0,R[$. Therefore, $u_\varepsilon$
is decreasing in $[0,R]$\,.\\
{\it Step~3.  The function $u_\varepsilon$ is decreasing in
$[R,1]$\,.}\\
Notice that from Eq. (\ref{u-Eq2}), any critical point $r_0\in]R,1[$
of $u_\varepsilon$ is a local minimum. Thanks to Steps~1~and~2, we
have also that $u'_\varepsilon(R)<0$ and $u''_\varepsilon(R)>0$.\\
Let us define the following function
$$f_\varepsilon(r)=\left\{
\begin{array}{lll}
\displaystyle\frac{u''_\varepsilon(R)}{2}\,(r-R)^2+u'_\varepsilon(R)\,(r-R)+u_\varepsilon(R)
&,&{\rm if~}0<r<R\,,\\
u_\varepsilon(r)&,&{\rm if~} R\leq r\leq 1\,,\\
f_\varepsilon(2-r)&,&{\rm if~} 1<r\leq 2\,.
\end{array}\right.$$
It is clear that $f_\varepsilon\in C^2([0,2])$ and that it satisfies
the following properties: (i) $r_0=1$ is a local minimum of
$f_\varepsilon$\,;
(ii)  if $r_0\in]0,2[$ is a critical point of $f_\varepsilon$, then
$r_0$ is a local minimum. This proves that $r_0=1$ is the only
critical point of $f_\varepsilon$ in $[0,2]$, and $f'_\varepsilon$
has a constant sign in $[0,1[$. Since $u'_\varepsilon(R)<0$, we
deduce that $u'_\varepsilon<0$ in $]R,1[$, hence the function
$u_\varepsilon$ is decreasing.
\hfill$\Box$\\

Finally, we mention without proof that the energy $C_0(\varepsilon)$ (cf.
(\ref{V-C0})) has the order of $\varepsilon^{-1}$, and we refer to
the methods in \cite[Section~6]{kach3} which permit to obtain the
leading order asymptotic expansion
$$C_0(\varepsilon)=\frac{c_1(a)}\varepsilon+c_2(a,R)+o(1),\quad(\varepsilon\to0),$$
where $c_1(a)$ and $c_2(a,R)$ are positive explicit constants.

\subsection{The case with magnetic field}\label{V-sec-minimizers}
This section is devoted to a preliminary analysis of the minimizers
of (\ref{V-EGL}) when $H\not=0$. The main point that we shall show
is how to extract the singular term $C_0(\varepsilon)$ (cf.
(\ref{V-C0})) from the energy of a minimizer.

Notice that the existence of minimizers is standard starting from a
minimizing sequence (cf. e.g. \cite{Gi}). A standard choice of gauge
permits one to assume that the magnetic potential satisfies
\begin{equation}\label{V-gauge}
{\rm div}\,A=0\quad {\rm in}~\Omega,\quad n\cdot A=0\quad{\rm
on}~\partial\Omega,
\end{equation}
where $n$ is the outward unit normal vector of
$\partial\Omega$.\\
With this choice of gauge, one is able to prove (since the boundaries
of $\Omega$ and $S_1$ are  smooth) that a minimizer $(\psi,A)$ is in $
C^1(\overline\Omega;\mathbb C)\times C^1(\overline\Omega;\mathbb
R^2)$. One  has also the following regularity (cf.
\cite[Appendix~A]{kach3}),
$$\psi \in C^2(\overline S_1;\mathbb C)\cup C^2(\overline S_2;\mathbb
C),\quad A\in C^2(\overline S_1;\mathbb R^2)\cup C^2(\overline
S_2;\mathbb R^2).$$

The next lemma is inspired from the work of Lassoued-Mironescu (cf.
\cite{LaMi}).

\begin{lem}\label{V-lem-psi<u}
Let $(\psi,A)$ be a minimizer of (\ref{V-EGL}). Then
$0\leq|\psi|\leq u_\varepsilon$ in $\Omega$, where $u_\varepsilon$
is the positive
minimizer of (\ref{V-EnH=0}).\\
Moreover, putting $\varphi=\frac{\psi}{u_\varepsilon}$, then the
energy functional (\ref{V-EGL}) splits in the form~:
\begin{equation}\label{V-splittingEn}
\mathcal G_{\varepsilon,H}(\psi,A)=C_0(\varepsilon)+ \mathcal
F_{\varepsilon,H}(\varphi,A),\end{equation} where $C_0(\varepsilon)$
has been introduced in (\ref{V-C0}) and the new functional $\mathcal
F_{\varepsilon,H}$ is defined by~:
\begin{equation}\label{V-reducedfunctional}
\mathcal F_{\varepsilon,H}(\varphi,A)=\int_\Omega
\left(u_\varepsilon^2
|(\nabla-iA)\varphi|^2
+\frac1{2\varepsilon^2}
u_\varepsilon^4(1-|\varphi|^2)^2+ |{\rm
curl}\,A-H|^2\right)\md x.\nonumber
\end{equation}
\end{lem}
\paragraph{\bf Proof}\ \\
The equality (\ref{V-reducedfunctional}) results from a direct but some
how long calculation, which permits to deduce in particular that
$\varphi$ is a solution of the equation
$$-(\nabla-iA)u_\varepsilon^2(\nabla-iA)\varphi=
\frac{u_\varepsilon^4}{\varepsilon^2}(1-|\varphi|^2)^2\varphi\,.$$\\
{\it Proof of $|\psi|\leq u_\varepsilon$.}\\
It is sufficient to prove that $|\varphi|\leq 1$. We shall invoke an
energy argument which we take from \cite{DGP}.\\
Let us introduce the set
$$\Omega_+=\{x\in\overline\Omega~:~|\varphi(x)|>1\}\,,$$
together with the functions (defined in $\Omega_+$)~:
$$f=\frac{\varphi}{|\varphi|}\,,\quad
\widetilde\varphi=[\,|\varphi|-1]_+f\,.$$
Then, it results from a direct calculation together with the
weak-formulation of the equation satisfied by $\varphi$ that
$$
0=\int_{\Omega_+}
\left(|\nabla|\varphi|\,|^2+(|\varphi|-1)
|\varphi|\,|(\nabla-iA)f|^2
+\frac1{2\varepsilon^2}
\left(1+|\varphi|)(1-|\varphi|)^2|\varphi|\right)u_\varepsilon^2\right)u_\varepsilon^2\,\md
x\,.
$$
Therefore, this yields that $\Omega_+$ has measure $0$.\hfill$\Box$\\

The estimate of the next lemma is very useful for exhibiting a
vortex-less regime for  minimizers  of (\ref{V-EGL}). It is due to
B\'ethuel-Rivi\`ere \cite{BeRi}, but a proof may be found also in
\cite[Corollary~3.1]{SaSe}
(see also \cite[Lemma~3.6]{kach4}).

\begin{lem}\label{grad-Abd}
Let $(\psi,A)$ be a minimizer of (\ref{V-EGL}). There exist
constants $C>0$ and $\varepsilon_0\in]0,1]$ such that, if the applied magnetic
field satisfies $H\ll\frac1{\varepsilon}$, then we have
$$|(\nabla-iA)\psi|\leq\frac{C}\varepsilon,\quad\forall~\varepsilon\in]0,\varepsilon_0]\,.$$
\end{lem}

Now, Lemma~\ref{grad-Abd} permits to conclude the following result.

\begin{lem}\label{BBH-thm3.3}
Assume that $(\psi,A)$ is a minimizer of (\ref{V-EGL}) and
let $\varphi=\frac\psi{u_\varepsilon}$. There exists a constant
$\mu_0>0$ such that if
$$\frac1{\varepsilon^2}\int_\Omega (1-|\varphi|^2)^2\,\md x\leq \mu_0\,,$$
then $|\varphi|\geq \frac12$ in $\overline\Omega$\,.
\end{lem}
\paragraph{\bf Proof.}
Lemma~\ref{grad-Abd} and the diamagnetic inequality together  yield
that
$$|\nabla |\psi|\,|\leq |(\nabla-iA)\psi|\leq\frac{C}\varepsilon,\quad
{\rm in}~\Omega\,.$$
Now, since
$$|\nabla u_\varepsilon|\leq \frac{C}\varepsilon\,$$
we deduce that
$$|\nabla|\varphi|\,|\leq \frac C{\varepsilon}\quad {\rm in}~\overline \Omega\,.$$
Thus, the result of the lemma becomes a consequence of
\cite[Theorem~III.3]{BBH}.\hfill$\Box$\\

\section{Analysis of the Meissner state}\label{V-Sec-meissnerstate}
Let us recall the definition of $u_\varepsilon$ and
$C_0(\varepsilon)$ in Theorem~\ref{V-thm-kach3} and (\ref{V-C0})
respectively. This section is devoted to the analysis of the
following variational problem (\ref{V-EM})~:
\begin{equation}\label{V-EM}
M_0(\varepsilon,H)=\min_{A\in H^1(\Omega;\mathbb R^2)} \mathcal
G_{\varepsilon,H}(u_\varepsilon,A),
\end{equation}
Since the function $u_\varepsilon$ is real-valued, one gets, for any
vector field $A$, the following decomposition~:
$$\mathcal G_{\varepsilon,H}(u_\varepsilon,A)=
C_0(\varepsilon)+\int_\Omega \left(|Au_\varepsilon|^2+ |{\rm
curl}\,A-H|^2\right)\,\md x.$$ Putting further
$$A=H\,\mathcal A,$$
\begin{equation}\label{V-J0}
J_0(\varepsilon)=\inf_{\mathcal A\in H^1(\Omega;\mathbb R^2)}\left[
\int_\Omega\left(|\mathcal A\,u_\varepsilon|^2+|{\rm curl}\,\mathcal
A-1|^2 \right)\,\md x\right],\end{equation} we get that
$$M_0(\varepsilon,H)=\inf_{A\in H^1(\Omega;\mathbb R^2)}
\mathcal G_{\varepsilon,H}(u_\varepsilon,A)=C_0(\varepsilon)+H^2J_0(\varepsilon),$$
and we are reduced to the analysis of the variational problem
(\ref{V-J0}).\\
Starting from a minimizing sequence (cf. \cite{SaSe}), it is
standard to prove that a minimizer $A_\varepsilon$ of (\ref{V-J0})
exists and satisfies the Coulomb gauge condition:
$${\rm div}\,A_\varepsilon=0\quad{\rm in}~\Omega,\quad n\cdot
 A_\varepsilon=0\quad{\rm on}~\partial\Omega,$$
where $n$ is the unit outward normal vector of the boundary of $\Omega$.\\
Notice also that $A_\varepsilon$ satisfies the Euler-Lagrange
equations~:
\begin{equation}\label{V-J0-EulerEq}
\nabla^\bot{\rm
curl\,}A_\varepsilon=u_\varepsilon^2\,A_\varepsilon\quad {\rm
in}~\Omega,\quad {\rm curl}\,A_\varepsilon=1\quad{\rm
on}~\partial\Omega.
\end{equation}
Here $\nabla^\bot=(-\partial_{x_2},\partial_{x_1})$
is the {\it Hodge gradient}.\\
Putting $h_\varepsilon={\rm curl}\,A_\varepsilon$, we get from the
first equation in (\ref{V-J0-EulerEq}) that
$A_\varepsilon=\frac1{u_\varepsilon^2}\nabla^\bot h_\varepsilon$. We
get also that $h_\varepsilon$ satisfies the equation:
\begin{equation}\label{V-hepsilon'}
-{\rm div}\left(\frac1{u_\varepsilon^2} \nabla
h_\varepsilon\right)+h_\varepsilon=0\quad{\rm in}~\Omega,\quad
h_\varepsilon=1\quad {\rm on}~\partial\Omega.
\end{equation}

\begin{lem}\label{V-lem-hepsilon}
The function $h_\varepsilon$ satisfies $0<h_\varepsilon<1$ in
$\Omega$, and it is the only function solving
(\ref{V-hepsilon'}).\\
Moreover, given $R'\in]0,R[$,
there exists a constant $c_0\in]0,1[$, and for each
$a\in\mathbb R_+\setminus\{1\}$, there exists a positive constant
$\varepsilon_0<1$   such that,
\begin{equation}\label{V-Eq-normhepsilon}
c_0\leq|h_\varepsilon(x)-1|< 1,\quad
\forall~x\in D(0,R'),~\forall~a\in\mathbb R_+\setminus\{1\},~
\forall~\varepsilon\in]0,\varepsilon_0]\,.
\end{equation}
\end{lem}
\paragraph{\bf Proof.}
The property that $0<h_\varepsilon<1$ and the uniqueness of
$h_\varepsilon$ are direct applications of the Strong Maximum
Principle.\\
Let us now prove (\ref{V-Eq-normhepsilon}).
Let us take a set $K\subset S_1$ (independent of
$\varepsilon$).
Due to the asymptotic behaviour of $u_\varepsilon$
(it remains exponentially close to $1$ in $K$, see
Lemma~\ref{interior}), it
results from (\ref{V-hepsilon'}) that $h_\varepsilon$ is bounded in
the $C^2$-norm of $K$. Thus, one can extract a subsequence of
$h_\varepsilon$, still denoted by $h_\varepsilon$, that converges to
a function $h\in C^2(K)$. The function $h$ satisfies the limiting
equation,
$$-\Delta h+h=0\quad{\rm in}~K.$$
By the Strong Maximum Principle, $0<h<1$ in $K$.
Let  $h_0$ be the solution of the equation
$$-\Delta h_0+h_0=0\quad{\rm in}~K,\quad h_0=1\quad{\rm on}~\partial
K.$$
Then, by the Strong Maximum Principle, $0<h\leq h_0<1$ in
$K$. This
 achieves the proof of the lemma.
\hfill$\Box$\\

\begin{lem}\label{V-lem-disc}
With the assumption (\ref{V-S}), the function
$h_\varepsilon$ is radial, i.e. $h_\varepsilon(x)=\widetilde
h_\varepsilon(|x|)$, with $\widetilde h_\varepsilon$ being an
increasing function.
\end{lem}
\paragraph{\bf Proof.}
That $h_\varepsilon$ is radial follows by the uniqueness of the
solution of (\ref{V-hepsilon'}) and by the fact that $u_\varepsilon$
is also radial.\\
The solution $h_\varepsilon$ being radial, i.e.
$$h_\varepsilon(x)=\widetilde h_\varepsilon(|x|),\quad
\forall~x\in\Omega,$$ let us  show  that the function $\widetilde
h_\varepsilon$ is increasing. For simplicity of notation, we shall
remove the tilde and write $h_\varepsilon$ for $\widetilde
h_\varepsilon$. Notice that $h_\varepsilon$ satisfies the
differential equation~:
\begin{equation}\label{V-Eq-radialh}
\left\{\begin{array}{l} - h''_\varepsilon(r)-\displaystyle\frac1{r}
h'_\varepsilon(r)
+2\frac{u'_\varepsilon(r)}{u_\varepsilon(r)}\,h_\varepsilon'(r)+u_\varepsilon^2(r)\,h_\varepsilon(r)=0,\quad
0<r<1,\\
h_\varepsilon'(0)=0,\quad h_\varepsilon(1)=1.
\end{array}\right.
\end{equation}
Let us calculate $h''_\varepsilon(0)$. Since $h'_\varepsilon(0)=0$,
we have  $h''_\varepsilon(0)=\displaystyle\lim_{r\to0}
\displaystyle\frac{h'_\varepsilon(r)}{r}$. Substituting in
(\ref{V-Eq-radialh}), we get that
\begin{equation}\label{V-h''(0)}
h_\varepsilon''(0)=\frac12u_\varepsilon^2(0)\,h_\varepsilon(0)>0.
\end{equation}
Let us introduce the even extension of $h_\varepsilon$, namely the
function
$$f_\varepsilon(r)=\left\{\begin{array}{l}
h_\varepsilon(r)\quad(r>0),\\
h_\varepsilon(-r)\quad (r<0).
\end{array}\right.$$
Then $f_\varepsilon$ satisfies the equation,
\begin{equation}\label{V-Eq-fepsilon}
-f''_\varepsilon(r)-\displaystyle\frac1{|r|} f'_\varepsilon(r)
+2\frac{\widetilde u'_\varepsilon(r)}{\widetilde
  u_\varepsilon(r)}\,f_\varepsilon'(r)+\widetilde
u_\varepsilon^2(r)\,f_\varepsilon(r)=0,\quad
r\in]-r_2,r_2[\setminus\{0\},
\end{equation}
and it attains a local minimum at $0$. We emphasize also here that
$\widetilde u_\varepsilon$ denotes the even extension of
$u_\varepsilon$.\\
If $r_0\in]-1,1[$ (with $r_0\not=0$) is a critical point of
$f_\varepsilon$, then it follows from (\ref{V-Eq-fepsilon}) that~:
$$f''_\varepsilon(r_0)=\widetilde u_\varepsilon^2(r_0)\,f_\varepsilon(r_0)>0.$$
If $r_0=0$, the conclusion $f''_\varepsilon(0)>0$ still holds,
thanks
to (\ref{V-h''(0)}).\\
Now these observations lead to the conclusion that $f_\varepsilon$
attains its minimum at a unique point, and that this point is the
only  critical point for $f_\varepsilon$. As we know that
$f'_\varepsilon(0)=0$, we get that $f_\varepsilon$ attains its
minimum at $0$ and that it is increasing
in $[0,1[$. This achieves the proof of the lemma.\hfill$\Box$\\

As we shall see, the next lemma
will play a distinguished role in the control of the
minimizing energy of `vortex balls'.

\begin{lem}\label{V-lem-controlgrad}
The following estimate holds
\begin{equation}\label{V-Eq-controlgrad}
\left\|\frac1{u_\varepsilon^2}\nabla
h_\varepsilon\right\|_{L^\infty(\Omega)}\leq 1, \quad
\forall~\varepsilon\in]0,1]\quad\forall~a\in\mathbb R_+\setminus\{1\}\,.
\end{equation}
\end{lem}
\subsubsection*{ Proof.} Notice that by Lemma~\ref{V-lem-disc},
$h_\varepsilon$ is radial. Then the equation for $h_\varepsilon$ can
be written in the form:
$$-\left(\frac{h'_\varepsilon}{u_\varepsilon^2}\right)'(r)
-\frac1r\,\frac{h_\varepsilon'}{u_\varepsilon^2}(r)
+h_\varepsilon(r)=0,\quad\forall~r\in]0,1[.
$$
Integrating this equation between $0$ and $r\in]0,1[$ and using the
fact that $h_\varepsilon$ is increasing, $h_\varepsilon'\geq0$, we
obtain:
$$\left(\frac{h'_\varepsilon}{u_\varepsilon^2}\right)(r)\leq
\int_{0}^rh_\varepsilon(\widetilde r)\,\md \widetilde r\leq
r\|h_\varepsilon\|_{L^\infty([0,1])}\leq 1,$$ which is the result of
the lemma.
\hfill$\Box$\\

Let us introduce the set
\begin{equation}\label{location-vortex}
\Lambda_\varepsilon=\left\{
x\in\Omega~:~\frac{1-h_\varepsilon(x)}{u_\varepsilon^2(x)}
=\max_{\overline\Omega}\frac{1-h_\varepsilon}{u_\varepsilon^2}\right\}\,.
\end{equation}

\begin{thm}\label{corollary}
The following two assertions hold.
\begin{enumerate}
\item If $a>1$, the function
  $\displaystyle\frac{1-h_\varepsilon}{u_\varepsilon^2}$
is strictly decreasing, and $\Lambda_\varepsilon=\{0\}$.
\item There exists $a_0\in]0,1[$ such that, for all $a\in]0,a_0[$, the
set $\Lambda_\varepsilon$ is a circle $\partial
D(0,R_\varepsilon)$
 localized strictly in $S_2$ as
$\varepsilon\to0$ in the sense that given $\alpha\in]0,1[$,
we have,
\begin{equation}\label{localization}
\varepsilon\ll R_\varepsilon-R\ll\varepsilon^\alpha\,,
\quad(\varepsilon\to0).
\end{equation}
Moreover, there exists a positive constant $C>0$ such that
\begin{equation}\label{control-u'}
|\nabla u_\varepsilon(x)|\leq C,\quad\forall~x\in S_2\setminus
 D(0,R_\varepsilon)\,.
\end{equation}{\rm }
\end{enumerate}
\end{thm}
\paragraph{\bf Proof.}
The proof of the first assertion is straightforward: When $a>1$,
the functions $u_\varepsilon$  and $h_\varepsilon$ are strictly
increasing, hence
$$\left(\frac{1-h_\varepsilon}{u_\varepsilon^2}  \right)'
=-\frac{u_\varepsilon
h'_\varepsilon+2(1-h_\varepsilon)u'_\varepsilon}
{u_\varepsilon^3}<0\,.$$ The proof of the second  assertion of the
corollary  is more
delicate. We present it in five steps.\\
{\it Step.~1. Proof of $\varepsilon\ll R_\varepsilon-R$.}\\
 Choose
$x_\varepsilon\in\Lambda_\varepsilon$ and let
$r_\varepsilon=|x_\varepsilon|$. Then $r_\varepsilon\in[0,1[$.
Thanks to Lemmas~\ref{V-lem-hepsilon}-\ref{V-lem-controlgrad}, we
have:
\begin{equation}\label{properties-max}
\liminf_{\varepsilon\to0}
\frac{1-h_\varepsilon(r_\varepsilon)}{u_\varepsilon^2(r_\varepsilon)}\geq
c_0,\quad
\limsup_{\varepsilon\to0} r_\varepsilon<1\,.
\end{equation}
Since $r_\varepsilon$ is a critical point of the function
$\frac{1-h_\varepsilon}{u_\varepsilon^2}$, we have
\begin{equation}\label{cp}
u'_\varepsilon(r_\varepsilon)=\frac{u_\varepsilon(r_\varepsilon)
h'_\varepsilon(r_\varepsilon)}{2(1-h_\varepsilon(r_\varepsilon))}\,.
\end{equation}
Then, by Lemma~\ref{V-lem-controlgrad},
$|u'_\varepsilon(r_\varepsilon)|\leq C$
for an explicit constant $C>0$. By Theorem~\ref{mainthm-H=0}, we
should have
$$|R-r_\varepsilon|\gg \varepsilon\quad{\rm as}~\varepsilon\to0\,.$$
Assume by
contradiction that $r_\varepsilon<R$. Then Theorem~\ref{mainthm-H=0}
yields that
$\displaystyle\lim_{\varepsilon\to0}u_\varepsilon(r_\varepsilon)=1$.\\
Let $\alpha\in]0,1[$ and choose $r'\in]0,R[$ such that
$$|h_\varepsilon(R+\varepsilon^\alpha)-h_\varepsilon(r')|\leq
\frac12|1-h_\varepsilon(r')|\,.$$
This choice of $r'$ is always possible, thanks to
Lemmas~\ref{V-lem-hepsilon}~and~\ref{V-lem-controlgrad}.\\
Now, notice that, as $\varepsilon\to0$,
\begin{eqnarray}\label{*100}
\frac{1-h_\varepsilon(R+\varepsilon^\alpha)}
{u_\varepsilon^2(R+\varepsilon^\alpha)}&\geq&
\frac12\,\frac{|1-h_\varepsilon(r')|}{u_\varepsilon^2(R+\varepsilon^\alpha)}
\nonumber\\
&\geq&\frac{c_0}{2\,u_\varepsilon^2(R+\varepsilon^\alpha)}
\quad[c_0\in]0,1[~{\rm given~in~Lemma~\ref{V-lem-hepsilon}]}\nonumber\\
&=&\frac{c_0}{2a}(1+o(1))\quad
\,{\rm [by~Theorem~\ref{mainthm-H=0}].}
\end{eqnarray}
On the other hand, by the definition of $r_\varepsilon$,
$$\frac{1-h_\varepsilon(R+\varepsilon^\alpha)}
{u_\varepsilon^2(R+\varepsilon^\alpha)}\leq
\frac{1-h_\varepsilon(r_\varepsilon)}{u_\varepsilon^2(r_\varepsilon)}\,,$$
and since
$\displaystyle\lim_{\varepsilon\to0}u_\varepsilon(r_\varepsilon)=1$,
we get
$$ \frac{1-h_\varepsilon(R+\varepsilon^\alpha)}
{u_\varepsilon^2(R+\varepsilon^\alpha)}\leq 1+o(1)\quad{\rm as~}
\varepsilon\to0\,.$$
Therefore, by choosing $a\in]0,\frac{c_0}2[$,
(\ref{*100}) leads to a contradiction. By putting
$$R_\varepsilon=\min_{x\in\Lambda_\varepsilon}|x|\,,$$
we
get the desired  statement: $\varepsilon\ll R_\varepsilon-R$ as
$\varepsilon\to0$.\\
{\it Step~2. Proof of $R_\varepsilon-R\ll 1$\,.}\\
Assume that there exists $R_1\in]R,1[$ such that, up to extraction
of a subsequence, $R_\varepsilon\to R_1$ as $\varepsilon\to0$. Let
$\delta=\frac12\min(|R-R_1|,1)$. We may assume, by extracting a
subsequence,
that $h_\varepsilon(R_1-\delta)\to c_*$ for some constant
$c_*\in]0,1[$. Then, by standard elliptic estimates, there exists a
    function
$h_*\in C^2\left(\overline{D(0,1)\setminus D(0,R_1-\delta)}\right)$
such that, upon the extraction of a
subsequence, we have,
$$h_\varepsilon\to h_*\quad{\rm in}~
C^2\left(\overline {D(0,1)\setminus D(0,R_1-\delta)}\right)\,,$$
and $h_*$ is a radial function and the unique solution of the equation
$$\left\{
\begin{array}{l}
-\Delta h_*+a\,h_*=0\quad{\rm in }~D(0,1)\setminus \overline{
D(0,R_1-\delta)}\,,\\
h_*=c_*\quad{\rm on}~\partial D(0,R_1-\delta),\quad h_*=1\quad{\rm
  on}~\partial D(0,1)\,.\end{array}\right.$$
A simple application of the maximum principle yields that $c_*<h_*<1$
in $D(0,1)\setminus \overline{
D(0,R_1-\delta)}$. Therefore, there exists a constant $M>0$ such that
$$h_*(R_1-\delta)<h_*(R_1)-M\,.$$
Consequently, when $\varepsilon$ is sufficiently small, we get the
lower bound:
$$\frac{1-h_\varepsilon(R_1-\delta)}{u_\varepsilon^2(R_1-\delta)}
>\frac{1-h_\varepsilon(R_1)}{u_\varepsilon^2(R_1)}+\frac{M}2\,,$$
and the same estimate holds when we replace $R_1$ by $R_\varepsilon$
and $\frac{M}2$ by $\frac{M}4$. This contradicts the definition
of $R_\varepsilon$, proving thus the desired property:
$R_\varepsilon-R\ll 1$ as $\varepsilon\to0$.\\
{\it Step~3. Finer localization:
Proof of $R_\varepsilon-R\ll\varepsilon^\alpha$\,.}\\
Assume that there exists $\alpha\in]0,1[$ such that, up to the
extraction of a subsequence, $R_\varepsilon>R+\varepsilon^\alpha$.\\
Let $\alpha'\in]\alpha,1[$ and set $\delta_\varepsilon=R_\varepsilon-R
-\varepsilon^{\alpha'}$. Notice that
$$\delta_\varepsilon\geq\frac{\varepsilon^{\alpha}}2\quad{\rm
  when~}\varepsilon
~{\rm is ~small~enough.}$$
Thanks to (\ref{cp}) and Lemma~\ref{interior},
  $h'_\varepsilon(R_\varepsilon)$
is exponentially small as $\varepsilon\to0$. Thus, from the equations
  satisfied by $h_\varepsilon$, we may assume that up to the
  extraction of a subsequence,
$$h''_\varepsilon(R_\varepsilon)\to \lambda_0~{\rm
    as~}\varepsilon\to0,
\quad\lambda_0>0\,.$$
Now, applying Taylor's formula to the function $h_\varepsilon$, we get
$$h_\varepsilon(R_\varepsilon-\delta_\varepsilon)=h_\varepsilon(R_\varepsilon)
+\lambda_0\delta_\varepsilon^2+o(\delta_\varepsilon^2)
\quad{\rm as}~\varepsilon\to0\,.$$
Consequently, thanks again to Lemma~\ref{interior}, we deduce that
\begin{eqnarray*}
\frac{1-h_\varepsilon(R_\varepsilon-\delta_\varepsilon)}
{u_\varepsilon^2(R_\varepsilon-\delta_\varepsilon)}&=&
\frac{1-h_\varepsilon(R_\varepsilon)}{u_\varepsilon^2(R_\varepsilon)}-\lambda_0\delta_\varepsilon^2+o(\delta_\varepsilon^2)\\
&<&\frac{1-h_\varepsilon(R_\varepsilon)}{u_\varepsilon^2(R_\varepsilon)}
-\frac{\lambda_0}2\delta_\varepsilon^2\,.
\end{eqnarray*}
Since the function $[0,1]\ni r\mapsto
\frac{1-h_\varepsilon(r)}{u_\varepsilon^2(r)}$ achieves its maximum on
$R_\varepsilon$, we get a contradiction. Therefore, we have proved the
desired localization for $R_\varepsilon$: Given $\alpha\in]0,1[$,
$R_\varepsilon-R\ll\varepsilon^\alpha$ as $\varepsilon\to0$.\\
{\it Step~4. Upper Bound for $|\nabla u_\varepsilon|$.}\\
Let us prove now (\ref{control-u'}). We know that
$|u'_\varepsilon(R_\varepsilon)|\leq C$, for some explicit constant $C>0$. On the
other hand,  by  Lemma~\ref{monotonicity}, the
function $u_\varepsilon$ is decreasing when $a<1$, hence
$u'_\varepsilon\leq 0$ in $]R,1]$. So it
is sufficient to prove that $u'_\varepsilon$ is increasing in
$[R,1]$.  Actually, coming back to the
equation of $u_\varepsilon$, we have, thanks to Theorem~\ref{V-thm-kach3},
$$u''_\varepsilon=-\frac1ru'_\varepsilon-\frac1{\varepsilon^2}
(a-u_\varepsilon^2)u_\varepsilon>0\quad{\rm in}~]R,1[\,,$$
hence we have the desired property regarding the monotonicity of
$u'_\varepsilon$. This achieves the proof of (\ref{control-u'}).\\
{\it Step~5. The function
$[0,1]\ni  r\mapsto\displaystyle\frac{1-h_\varepsilon(r)}{u_\varepsilon^2(r)}$
  achieves its maximum on a unique point.}\\
Let us prove now that $\Lambda_\varepsilon=\partial
D(0,R_\varepsilon)$, i.e. the radial function
$\frac{1-h_\varepsilon}{u_\varepsilon^2}$ attains its  maximum
uniquely at $R_\varepsilon$.\\
By Lemma~\ref{V-lem-hepsilon}, there exists a constant $R_*\in]R,1[$
such that any  maximum point $x\in\Lambda_\varepsilon$ should satisfy $R<|x|<R_*$.  Let $r_\varepsilon\in]R_\varepsilon,R_*[$ be a critical point of
$\displaystyle
\frac{1-h_\varepsilon}{u_\varepsilon^2}$. Then,
$$\left(\frac{1-h_\varepsilon}{u_\varepsilon^2}\right)''
=\frac{-3f_\varepsilon}{u_\varepsilon^4}\,,$$
where
$$f_\varepsilon=u'_\varepsilon h'_\varepsilon-\frac1rh'_\varepsilon+u_\varepsilon^2h_\varepsilon+2(1-h_\varepsilon)u''_\varepsilon\,.$$
It is sufficient to  prove that $f_\varepsilon(r_\varepsilon)>0$. We
distinguish between two cases:
$$({\rm i})~\limsup_{\varepsilon\to0}
u''_\varepsilon(r_\varepsilon)=\infty,\quad{\rm or}\quad
({\rm ii})~\limsup_{\varepsilon\to0}u''_\varepsilon(r_\varepsilon)<\infty\,.$$
In case (i), since $u'_\varepsilon$ is bounded in
$[R_\varepsilon,1[$,
we deduce easily that as $\varepsilon\to0$,
$$f_\varepsilon(r_\varepsilon)>0\,.$$
In case (ii), it is easy to verify that $u_\varepsilon''$ is
decreasing in $]R_\varepsilon,1[$. Hence there exists a constant $C>0$
such that, up to the extraction of a
subsequence, $u''_\varepsilon(r)\leq C$ in $[r_\varepsilon,1[$.\\
By the mean value theorem, we deduce that
$$|u'_\varepsilon(r_\varepsilon)
-u'_\varepsilon(r_\varepsilon+\varepsilon^\alpha)|\leq
C\varepsilon^\alpha.$$
Thus, we get by Lemma~\ref{interior} that
$|u'_\varepsilon(r_\varepsilon)|\ll1$
as $\varepsilon\to0$, and consequently, we get by (\ref{cp}) that
$h_\varepsilon'(r_\varepsilon)\ll 1$ as $\varepsilon\to0$.\\
Now, this yields in this case that
$f_\varepsilon(r_\varepsilon)>0$. Therefore, we have proved all the
statements of the theorem.
\hfill$\Box$\\

Let us introduce the function
\begin{equation}\label{xi-epsilon}
\Omega\ni x\mapsto
\xi_\varepsilon(x)=\frac{h_\varepsilon(x)-1}{u_\varepsilon^2(x)}\,,
\end{equation}
together with
\begin{equation}\label{lambda-epsilon}
\lambda_\varepsilon=\max_{x\in\overline\Omega}|\xi_\varepsilon(x)|\,.
\end{equation}

\begin{corol}\label{finite-points}
Let $a_0\in]0,1[$ be the constant of
Theorem~\ref{corollary}. For all $a\in]0,a_0[$, there exist  positive
constants  $\delta_*$ and $\varepsilon_0$
 such that~:
$$\left\{
\begin{array}{l}
\forall~x\in\Omega\quad{\rm s.t.}~~
|\,|x|-R|^2\geq|\ln\varepsilon|^{-1/2},\quad\forall~\varepsilon\in]0,\varepsilon_0],
\\
\xi_\varepsilon(x)\geq
-\lambda_\varepsilon+\delta_*|\ln\varepsilon|^{-1/2}\,.
\end{array}\right.$$
\end{corol}
\paragraph{\bf Proof.}
We make the following claim:
\begin{equation}\label{claim1}
\exists~c_0>0\,,\quad
\xi_\varepsilon(x)\geq -\lambda_\varepsilon+c_0\quad{\rm when~}|x|\leq
R-|\ln\varepsilon|^{-1/4}\,.\end{equation}
Once we prove (\ref{claim1}), we deduce the conclusion of the corollary when
$|x|<R-|\ln\varepsilon|^{-1/4}$.\\
The proof of (\ref{claim1}) is rather easy. First, notice that,
putting $r_\varepsilon=R-\varepsilon^\alpha$ with $\alpha\in]0,1[$, we
    have by Lemma~\ref{V-lem-controlgrad}
$$h_\varepsilon(r_\varepsilon)\geq
    h_\varepsilon(R_\varepsilon)+\mathcal O (\varepsilon^\alpha)\,.$$
On the other hand, Theorem~\ref{V-lem-hepsilon} yields that the
    function
$r\mapsto
\displaystyle\frac{h_\varepsilon(r)-1}{u_\varepsilon^2(r)}$ is decreasing
in $[0,R]$. Thus, for all $r\in[0,R-|\ln\varepsilon|^{-1/4}]$,
we have
$$\frac{h_\varepsilon(r)-1}{u_\varepsilon^2(r)}\geq
 \frac{h_\varepsilon(r_\varepsilon)-1}{u_\varepsilon^2(r_\varepsilon)}.$$
Therefore,
$$\xi_\varepsilon(r)\geq
    \frac{h_\varepsilon(R_\varepsilon)-1}{u_\varepsilon^2(r_\varepsilon)}
+\mathcal O(\varepsilon^\alpha)\,.$$
Invoking Lemma~\ref{interior}, we deduce that
$$\xi_\varepsilon(r)\geq -\lambda_\varepsilon+c_0\,,$$
where
$c_0\in ]0,\left(\frac1a-1\right)
\displaystyle\liminf_{\varepsilon\to0}(1-h_\varepsilon(R_\varepsilon))[$.\\
Now, let us prove the conclusion  of the corollary when
$R+|\ln\varepsilon|^{-1/4}\leq|x|\leq1$.
By Lemma~\ref{interior},
it is sufficient to find $\delta>0$ and $r_0>0$ such that
\begin{equation}\label{claim2}
h_\varepsilon(r)\geq h_\varepsilon(R_\varepsilon)+\delta\max(
|r-R|^2,|\ln\varepsilon|^{-1/2})\,,\quad
R+|\ln\varepsilon|^{-1/4}\leq r\leq R+r_0\,.\end{equation}
To prove (\ref{claim2}), we deal separately with the case whether
$\displaystyle\liminf_{\varepsilon\to0}h'_\varepsilon(R_\varepsilon)=0$
or
$\displaystyle\liminf_{\varepsilon\to0}h'_\varepsilon(R_\varepsilon)>0$.\\
{\it Proof of (\ref{claim2}) when
  $\displaystyle\liminf_{\varepsilon\to0}h'_\varepsilon(R_\varepsilon)=0$.}\\
In this case, there exists $c_0>0$ such that, up to the extraction of
a subsequence,
$$h'_\varepsilon(R_\varepsilon)\to 0,\quad
h''_\varepsilon(R_\varepsilon)\to c_0\quad{\rm
  as}~\varepsilon\to0\,.$$
Set $r_\varepsilon=R+\varepsilon^\alpha$ where $\alpha\in]0,1[$ is
  given. By Theorem~\ref{corollary}, we have
$$h'_\varepsilon(r_\varepsilon)\to 0,\quad
h''_\varepsilon(r_\varepsilon)\to c_0\quad{\rm
  as}~\varepsilon\to0\,.$$
Moreover, by Lemma~\ref{interior} and the equation of $h_\varepsilon$,
$h'''_\varepsilon(r)$ is bounded in
$[r_\varepsilon,1]$.
Therefore, applying Taylor's formula
up to the order $2$, we get a positive constant $r_0\in]0,1[$ such that
\begin{eqnarray}\label{Taylor-hepsilon}
h_\varepsilon(r)&=&h_\varepsilon(r_\varepsilon)
+h'_\varepsilon(r_\varepsilon) (r-r_\varepsilon)
+\frac12h''_\varepsilon(r_\varepsilon)(r-r_\varepsilon)^2+\mathcal
O(|r-r_\varepsilon|^3)\nonumber\\
&\geq&
h_\varepsilon(r_\varepsilon)
+\frac{c_0}2(r-r_\varepsilon)^2\end{eqnarray}
provided that $0<r-r_\varepsilon<r_0$.\\
Thanks to Theorem~\ref{corollary},
$\varepsilon<R_\varepsilon-R<\varepsilon^\alpha$. Hence by
Lemma~\ref{V-lem-controlgrad},
$$h_\varepsilon(r_\varepsilon)=h_\varepsilon(R_\varepsilon)
+\mathcal O(\varepsilon^\alpha)\,.$$
Therefore, when $|\ln\varepsilon|^{-1/2}<r-R_\varepsilon<r_0$,
(\ref{Taylor-hepsilon}) is
nothing but
(\ref{claim2}).\\
{\it Proof of (\ref{claim2}) when
  $\displaystyle\liminf_{\varepsilon\to0}h'_\varepsilon(R_\varepsilon)>0$.}\\
We may assume in this case that
$h'_\varepsilon(R_\varepsilon)\to c_0>0$ as $\varepsilon\to0$.
By Theorem~\ref{corollary} and the equation of $h_\varepsilon$,
$h''_\varepsilon(r)$ is bounded in $[R_\varepsilon,1]$
independently of $\varepsilon$.\\
We apply again Taylor's formula
\begin{eqnarray*}
h_\varepsilon(r)&=&h_\varepsilon(R_\varepsilon)
+h'_\varepsilon(R_\varepsilon)
(r-R_\varepsilon)
+\mathcal
O(|r-R_\varepsilon|^2)\\
&\geq&
h_\varepsilon(R_\varepsilon)
+\frac{c_0}2|r-R_\varepsilon|\,,\end{eqnarray*}
This is nothing but again
(\ref{claim2}). This achieves the proof of the corollary.\hfill$\Box$\\

\section{Lower bound of the energy}\label{Section:LB}
\subsection{Vortex-balls}
In this section we construct suitable `vortex-balls' providing a lower
bound of the energy of minimizers of (\ref{V-EGL}). Recall the
decomposition of the energy in Lemma~\ref{V-lem-psi<u}, which permits
us to work with the `reduced energy functional' $\mathcal
F_{\varepsilon,H}$.\\
Notice that, by using
$\left(u_\varepsilon,\frac1{u_\varepsilon^2}\nabla^\bot
h_\varepsilon\right)$ as a test configuration for the function
(\ref{V-EGL}), we deduce the upper bound~:
\begin{equation}\label{V-upperbound}
\mathcal F_{\varepsilon,H}(\varphi,A)\leq H^2J_0(\varepsilon)\,,
\end{equation}
where $\varphi=\psi/u_\varepsilon$, $(\psi,A)$ a minimizer of
(\ref{V-EGL}),
and $J_0(\varepsilon)$ is introduced in (\ref{V-J0}),
\begin{equation}\label{V-J0'}
J_0(\varepsilon)=\int_\Omega\left(\frac1{u_\varepsilon^2}|\nabla
h_\varepsilon|^2+|h_\varepsilon-1|^2\right)\,\md x\,.
\end{equation}
We shall always work under the hypothesis that there exists a positive
constant $C>0$ such that the applied magnetic field $H$ satisfies
\begin{equation}\label{hypothesis-H}
H\leq C|\ln\varepsilon|\,.\end{equation}

The upper bound (\ref{V-upperbound}) provides us, as in \cite{SaSe}, with the
construction of suitable `vortex-balls'.

\begin{prop}\label{V-lem-vortexballs}
Let $(\psi,A)$ be a minimizer of (\ref{V-EGL}) and
$\varphi=\displaystyle\frac{\psi}{u_\varepsilon}$. Then, under the
hypotheses (\ref{hypothesis-H}), for each $p\in]1,2[$,
there exist a constant $C>0$
and a finite family of disjoint balls
$\{B((a_i,r_i)\}_{i\in I}$ satisfying the following properties:
\begin{enumerate}
\item $w=\{x\in \overline\Omega~:~|\varphi(x)|\leq
  1-|\ln\varepsilon|^{-4}\}
\subset\displaystyle\displaystyle\bigcup_{i\in I}B(a_i,r_i)$.
\item $\displaystyle\sum_{i\in I}r_i\leq
C\,|\ln\varepsilon|^{-10}$.
\item Letting $d_i$ be the degree of the
function $\varphi/|\varphi|$ restricted to $\partial B(a_i,r_i)$ if
$B(a_i,r_i)\subset\Omega$ and $d_i=0$
  otherwise, then  we have:
\begin{eqnarray}\label{V-Eq-LBest}
&&\hskip-0.5cm \int_{B(a_i,r_i)\setminus\omega}
u_\varepsilon^2|(\nabla-iA)\varphi|^2\,\md x+\int_{B(a_i,r_i)}|{\rm
  curl}\,A-H|^2\,\md x\geq\\
&&\hskip3.5cm 2 \pi|d_i|
\left(\min_{B(a_i,r_i)}u_\varepsilon^2\right)
\left(|\ln\varepsilon|-C\ln|\ln\varepsilon|\right).
\nonumber
\end{eqnarray}
\item
$\left\|2\pi\displaystyle\sum_{i\in I} d_i\delta_{a_i}-{\rm
    curl}\big{(}A+(i\varphi,\nabla_A\varphi)\big{)}
\right\|_{W^{-1,p}_0(\Omega)}\leq
\max(|\ln\varepsilon|^{-4}).$
\end{enumerate}
\end{prop}

We follow the usual terminology and call the balls constructed in
Proposition~\ref{V-lem-vortexballs} `vortex-balls'. The proof of Proposition~\ref{V-lem-vortexballs} is very similar to
that of \cite[Proposition~5.2]{kach4}, and is actually a simple
consequence of the analysis of \cite{SaSe}.\\

Proposition~\ref{V-lem-vortexballs} permits us to prove the following
theorem.

\begin{thm}\label{LowerBound-thm}
Let $(\psi,A)$ be a minimizer of (\ref{V-EGL}) and
$\varphi=\displaystyle\frac{\psi}{u_\varepsilon}$. Then, under the
hypothesis (\ref{hypothesis-H}),
there exist a constant $C>0$ and a finite family of disjoint balls
$\{B((a_i,r_i)\}_{i\in I}$ such that~:
\begin{enumerate}
\item $\displaystyle\sum_{i\in I} r_i\leq C|\ln\varepsilon|^{-10}$;\\
\item $|\varphi|\geq\frac12$ on $\Omega\setminus \cup_iB(a_i,r_i)$.
\item Letting $d_i$ be the degree of the
function $\varphi/|\varphi|$ restricted to $\partial B(a_i,r_i)$ if
$B(a_i,r_i)\subset\Omega$ and $d_i=0$
  otherwise, then  we have:
\begin{equation}\label{V-Eq-LBest'}
\begin{split}
\mathcal F_{\varepsilon,H}(\varphi,A)\,\geq& \,H^2J_0(\varepsilon)\\
&+\int_{\Omega\setminus\cup_i
  B(a_i,r_i)}\frac1{u_\varepsilon^2}|\nabla(h-H\,h_\varepsilon)|^2\,\md x
+\int_\Omega
|h-H\,h_\varepsilon|^2\,\md x\\
&+2\pi\sum_{i\in I}
\left[\left(\min_{B(a_i,r_i)}u_\varepsilon^2\right)(|\ln\varepsilon|
-C\ln|\ln\varepsilon|)\right]|d_i|\\
&+4\pi H\sum_{i\in I}d_i(h_\varepsilon-1)(a_i)-CH|\ln\varepsilon|^{-4}\,,\\
\end{split}
\end{equation}
where $h={\rm curl}\,A$ and $h_\varepsilon$ is introduced in
(\ref{V-hepsilon'}).
\end{enumerate}
\end{thm}

The proof is essentially  that of \cite[Theorem~5.3]{kach4}.\\

Let us recall the definition of $\lambda_\varepsilon$ in
(\ref{lambda-epsilon}). We put
\begin{equation}\label{k-epsilon}
k_\varepsilon=\frac1{2\lambda_\varepsilon}=
\frac12\left(\max_{x\in\overline\Omega}|\xi_\varepsilon(x)|\right)^{-1}\,.
\end{equation}

\begin{corol}\label{LowerBound-corol}
With the notations of Theorem~\ref{LowerBound-thm}, the following
lower bound holds:
\begin{equation}\label{V-Eq-LBest''}
\begin{split}
\mathcal F_{\varepsilon,H}(\varphi,A)\,\geq& \,H^2J_0(\varepsilon)\\
&+\int_{\Omega\setminus\cup_i
  B(a_i,r_i)}\frac1{u_\varepsilon^2}|\nabla(h-H\,h_\varepsilon)|^2\,\md x
+\int_\Omega
|h-H\,h_\varepsilon|^2\,\md x\\
&+2\pi\chi_\varepsilon(a)\sum_{d_i>0}
\left(|\ln\varepsilon|-2k_\varepsilon^{-1}H-C\ln|\ln\varepsilon|\right)d_i\\
&+\min(1,a)\sum_{d_i\leq0}(|\ln\varepsilon|-C\ln|\ln\varepsilon|)|d_i|-CH|\ln\varepsilon|^{-4}\,,\\
\end{split}
\end{equation}
where $\chi_\varepsilon(a)=\min(1,a)$ if
$|\ln\varepsilon|-2k_\varepsilon^{-1}H\geq0$, and
$\chi_\varepsilon(a)=\max(1,a)$ otherwise.
\end{corol}
\paragraph{\bf Proof.}
Let us assign to each ball $B(a_i,r_i)$ a point $a'_i\in
\overline{B(a_i,r_i)\cap\Omega}$ such that
$$u_\varepsilon(a'_i)=\min_{\overline{B(a_i,r_i)}}u_\varepsilon\,.$$
Then, thanks to Lemma~\ref{V-lem-controlgrad} and to the first point
of Theorem~\ref{LowerBound-thm}, there exists a constant $c>0$ such
that
$$\forall~i,\quad |h_\varepsilon(a_i)-h_\varepsilon(a_i')|\leq
c|a_i-a'_i|\leq c|\ln\varepsilon|^{-10}\,.$$ This permits us to
write
\begin{eqnarray*}
&&\hskip-0.5cm\sum_{d_i>0}\left[
\left(\min_{B(a_i,r_i)}u_\varepsilon^2\right)|\ln\varepsilon|+2H(h_\varepsilon(a_i)-1)\right]d_i\\
&&\geq \sum_{d_i>0}u_\varepsilon^2(a'_i)\left[
|\ln\varepsilon|-\left(\frac{1-h_\varepsilon(a'_\varepsilon)}{2u_\varepsilon^2(a'_i)}\right)H
-2c|\ln\varepsilon|^{-4}H\right]d_i\,.
\end{eqnarray*}
By definition of $k_\varepsilon$, we have
$\displaystyle\frac{1-h_\varepsilon(a'_\varepsilon)}{2u_\varepsilon^2(a'_i)}\leq
k_\varepsilon^{-1}$. By Theorem~\ref{V-thm-kach3}, $\min(1,a)\leq
u_\varepsilon^2(a_i')\leq\max(1,a)$. Therefore, we get
\begin{equation}\label{Proof-Corol:LB1}
\sum_{d_i>0}\left[
\left(\min_{B(a_i,r_i)}u_\varepsilon^2\right)|\ln\varepsilon|+2H(h_\varepsilon(a_i)-1)\right]d_i
\geq
\chi_\varepsilon(a)\left(|\ln\varepsilon|-2k_\varepsilon^{-1}H\right)d_i\,.
\end{equation}
For the terms with negative degrees, we write
\begin{equation}\label{Proof-Corol:LB2}
 \sum_{d_i\leq0}\left[
\left(\min_{B(a_i,r_i)}u_\varepsilon^2\right)|\ln\varepsilon|-2H(h_\varepsilon(a_i)-1)\right]|d_i|
\geq
\min(1,a)\sum_{d_i\leq 0}|\ln\varepsilon|\,|d_i|\,.
\end{equation}
Substituting (\ref{Proof-Corol:LB1})-(\ref{Proof-Corol:LB2}) in
(\ref{V-Eq-LBest'}), we get the desired lower bound of the
corollary.\hfill$\Box$\\

\subsection{Upper bound on the total degree}
Let us assume from now on that $(\psi,A)$ is a minimizer of
(\ref{V-EGL}) and that $(B(a_i,r_i))_i$ is its associated family of
vortex-balls constructed in Theorem~\ref{LowerBound-thm}.
Our aim is to
give an upper bound on the total degree $\sum_i|d_i|$. The answer will
be strongly dependent on the parameter $a$.\\

\begin{lem}\label{upperbound-degree}
Assume that for a given constant $K>0$, the magnetic field satisfies
$H\leq k_\varepsilon|\ln\varepsilon|+K\ln|\ln\varepsilon|$.
With  the notation of Theorem~\ref{LowerBound-thm},
the following two assertions hold.
\begin{enumerate}
\item If $a>1$, then there exists constants $C>0$ and
  $\varepsilon_0\in]0,1]$  such that,
\begin{equation}\label{a>1}
\max_{|d_i|>0}|a_i|\leq C|\ln\varepsilon|^{-1/4},\quad
\sum_i|d_i|\leq C,\quad \forall~\varepsilon\in]0,\varepsilon_0]\,.
\end{equation}
\item There exists $a_0\in]0,1[$ such that, if $a\in]0,a_0]$, there
exists positive constants $\varepsilon_0$ and $C$ such that,
\begin{equation}\label{a<1}
\max_{|d_i|>0} |R-|a_i|\,|\leq
C|\ln\varepsilon|^{-1/4},\quad\sum_i|d_i|\leq
C\ln|\ln\varepsilon|,\quad\forall~\varepsilon\in]0,\varepsilon_0]\,.
\end{equation}
\end{enumerate}
\end{lem}

Since the proof of Assertion (1) is very close to that of
\cite[Theorem~2]{SaSe4} (with  only very few technical modifications),
we omit it. We give rather the proof of the second assertion of the
lemma.\\
Let us introduce:
\begin{equation}\label{totaldegree:+-}
D_+=\sum_{i,\,d_i>0}|d_i|,\quad D_-=\sum_{i,\,d_i\leq0}|d_i|,\quad
D=D_++D_-=\sum_i|d_i|\,,
\end{equation}
and
\begin{equation}\label{totaldegree:D0}
D_0=\sum_{
|R-|a_i|\,|\leq|\ln\varepsilon|^{-1/4}
}|d_i|\,.
\end{equation}
We make the following two claims:
\begin{equation}\label{claim:D+}
\exists~C>0,\quad D_-\leq C\, D_+\frac{\ln|\ln\varepsilon|}
{|\ln\varepsilon|}\,,
\end{equation}
and
\begin{equation}\label{claim:D0}
\exists~C>0,\quad D-D_0\leq C\,D\frac{\ln|\ln\varepsilon|}{\sqrt{
|\ln\varepsilon|}}\,.
\end{equation}
Now we show that when the claims (\ref{claim:D+}) and (\ref{claim:D0})
hold, then we can prove Assertion (2) of
Lemma~\ref{upperbound-degree}.\\
We put $\widetilde \Omega=\Omega\setminus \cup_i B(a_i,r_i)$, where
$B(a_i,r_i)$ are the vortex-balls constructed in
Theorem~\ref{LowerBound-thm}.
For a given $t>0$, we denote by $C_t$ the circle of center $0$ and radius $t$,
and by $B_t$ the open ball of center $0$ and radius $t$.
Let us introduce the set of positive real numbers:
\begin{equation}\label{set:E}
E=\{t\in\,]R+|\ln\varepsilon|^{-1/4},1[~:~C_t\subset \widetilde\Omega\}\,.
\end{equation}
Thanks to Theorem~\ref{LowerBound-thm}, the set $E$ is non empty and
has a positive measure
$$\liminf_{\varepsilon\to0}|E|>0\,.$$
Theorem~\ref{LowerBound-thm} gives
 $|\varphi|\geq 1-|\ln\varepsilon|^{-4}$ on $C_t$ whenever $t\in
E$. Therefore, the degree
$$ d_t={\rm deg}\left(\frac{\varphi}{|\varphi|}, C_t\right)$$
is well defined whenever $t\in E$.\\
Writing $h={\rm curl}\,A$ and $\varphi=|\varphi|e^{i\phi}$
for an $H^2$-function $\phi$, the following equation holds
\begin{equation}\label{equation:h}
-\frac1{u_\varepsilon^2}\nabla^\bot
 h=|\varphi|^2(\nabla\phi-A)\quad{\rm in}~\widetilde\Omega\,.
\end{equation}
Let us recall also the equation for $h_\varepsilon$,
$$
-{\rm div}\left(\frac1{u_\varepsilon^2}\nabla
 h_\varepsilon\right)+h_\varepsilon=0\quad {\rm in}~\Omega\,,
$$
from which it follows, by Stoke's formula:
\begin{equation}\label{stoke}
\int_{C_t}\frac{1}{u_\varepsilon^2}\,n\cdot\nabla h_\varepsilon\,\md
\theta=\int_{B_t}h_\varepsilon\,\md x\,,
\end{equation}
where $n$ is the unit outward normal vector of $B_t$, $n(x)=x/|x|$ for
all $x\in\mathbb R^2\setminus\{0\}$.\\
On the other hand, it results from (\ref{equation:h}) and Stoke's formula,
\begin{eqnarray*}
\int_{C_t}\frac1{u_\varepsilon^2}\,n\cdot\nabla h\,\md\theta
&=&\int_{B_t}|\varphi|^2\,\tau\cdot(\nabla\phi-A)\,\md x\\
&=&\int_{B_t}\tau\cdot\nabla\phi\,\md x-\int_{C_t} h\,\md x+T(t)\,,
\end{eqnarray*}
where $(n,\tau)$ is a direct frame, and
\begin{equation}\label{T(t)}
T(t)=\int_{B_t}(|\varphi|^2-1)\,\tau\cdot(\nabla\phi-A)\,\md x\,.
\end{equation}
Coming back to the definition of the degree, we deduce that
\begin{equation}\label{Stoke:h}
\int_{C_t}\frac1{u_\varepsilon^2}\,n\cdot\nabla h\,\md\theta=
2\pi  d_t-\int_{B_t}h\,\md x+T(t)\,.
\end{equation}
Combining (\ref{stoke}) and (\ref{Stoke:h}), we get
$$\int_{C_t}\frac1{u_\varepsilon^2}\left(\nabla
h-H\nabla h_\varepsilon\right)\cdot n\,\md \theta+\int_{B_t}
(
h-H\,h_\varepsilon)\,\md x=2\pi d_t+T(t)\,.$$
Applying Cauchy-Schwarz inequality on each integral and squaring, we
obtain (recall that the function $u_\varepsilon$ is radial)
$$\int_{C_t}\frac1{u_\varepsilon^2}\left|\nabla(
h-H\, h_\varepsilon)\right|^2\,\md \theta+\frac{t}2\int_{B_t}
\left|
h-H\,h_\varepsilon\right|^2\,\md x\geq \frac{\pi}{3t} d_t^2-C\left[
T^2(t)+u_\varepsilon^{-2}(t)\right].$$
Thanks to (\ref{claim:D+}) and (\ref{claim:D0}), we infer from the
above lower bound
\begin{equation}\label{d_t=D}
\int_{C_t}\frac1{u_\varepsilon^2}\left|\nabla(
h-H\, h_\varepsilon)\right|^2\,\md \theta+\frac{t}2\int_{B_t}
\left|
h-H\,h_\varepsilon\right|^2\,\md x\geq \frac{\pi}{4t} D^2-C\left[
T^2(t)+u_\varepsilon^{-2}(t)\right]\,,\end{equation}
where $D$ is the total degree introduced in (\ref{totaldegree:+-}). \\
Now, we integrate both sides of (\ref{d_t=D}) with respect to $t$ and
we recall that $\inf E>R$. This yields
\begin{eqnarray*}
&&\hskip-1cm\int_{\widetilde \Omega}
\frac1{u_\varepsilon^2}\left|\nabla(
h-H\, h_\varepsilon)\right|^2\,\md x+\int_{\Omega}
\left|
h-H\,h_\varepsilon\right|^2\,\md x\\
&&\geq
\int_{E}\left(\int_{C_t}\frac1{u_\varepsilon^2}\left|\nabla(
h-H\, h_\varepsilon)\right|^2\,\md \theta+\frac{t}2\int_{B_t}
\left|
h-H\,h_\varepsilon\right|^2\,\md x\right)\md t\\
&&\geq\int_E\left(\frac{\pi}{4t} D^2-C\left[
T^2(t)+u_\varepsilon^{-2}(t)\right]\right)\,\md t\\
&&\geq \frac{|E|}{4}D^2-C
 \int_E\left(T^2(t)+u_\varepsilon^{-2}(t)\right)\,\md t\\
&&\geq \widetilde C\left[D^2- \int_E\left(T^2(t)+u_\varepsilon^{-2}(t)\right)\,\md t\right]\,,
\end{eqnarray*}
where $\widetilde C>0$ is an explicit constant.\\
Since $u_\varepsilon^{2}>a$ when $a<1$  (see
Theorem~\ref{V-thm-kach3}), it is clear that
$\int_Eu_\varepsilon^{-2}\,\md t\leq
 a^{-1}|E|\leq
C$. Let us estimate the integral of $T^2(t)$. Notice that
$$\int_E T^2(t)\,\md t\leq \int_{\widetilde\Omega}
(1-|\varphi|^2)|(\nabla-iA)\varphi|^2\,\md
x\leq |\ln\varepsilon|^{-4}\int_\Omega|(\nabla-iA)\varphi|^2\,\md
x\ll1$$
where we have used Theorem~\ref{LowerBound-thm} and the constraint on
the applied magnetic field $H=\mathcal O(|\ln\varepsilon|)$.\\
Therefore, we conclude finally that, for a possibly larger explicit
constant $\widetilde C>0$,
\begin{equation}\label{conclusion}
\int_{\widetilde \Omega}
\frac1{u_\varepsilon^2}\left|\nabla(
h-H\, h_\varepsilon)\right|^2\,\md x+\int_{\Omega}
\left|
h-H\,h_\varepsilon\right|^2\,\md x\geq \widetilde C(D^2-1)\,.
\end{equation}
We substitute (\ref{claim:D+}), (\ref{claim:D0}) and
(\ref{conclusion})
in (\ref{V-Eq-LBest'}) to obtain:
\begin{equation}\label{LowerBound-conclusion}
\mathcal F_{\varepsilon,H}(\varphi,A)\geq H^2J_0(\varepsilon)
+\widetilde C (D^2-1)-C\ln|\ln\varepsilon|D\,.
\end{equation}
Matching this lower bound with the upper bound (\ref{V-upperbound}), we deduce that
$$D^2\leq C'\ln|\ln\varepsilon|D\,,$$
which permits us to deduce the statement concerning the total degree
in the second assertion  of
Lemma~\ref{upperbound-degree}. Substituting the bound of $D$ in
(\ref{claim:D+}) and (\ref{claim:D0}), we deduce that
$$D_-=0,\quad D_0=D\,,$$
thus proving that all  the vortices have positive degrees together
with the first statement in the assertion (2) of
Lemma~\ref{upperbound-degree}.\\
We have only to prove Claims (\ref{claim:D+}) and (\ref{claim:D0}).
Claim (\ref{claim:D+}) is a direct consequence of
Theorem~\ref{LowerBound-thm}. Claim (\ref{claim:D0}) is a
consequence of Lemma~\ref{finite-points}.

\section{Upper bound of the energy}\label{Section:UB}
\subsection{Main result}
In this section, we assume that the magnetic field satisfies
\begin{equation}\label{H=Hc}
H=k_\varepsilon|\ln\varepsilon|+\lambda\ln|\ln\varepsilon|,\quad(\lambda\in\mathbb
R)\,,
\end{equation}
where $k_\varepsilon$ is introduced in (\ref{k-epsilon}).\\

The aim of this section is to establish the following upper bound for the
energy $\mathcal F_{\varepsilon,H}(\varphi,A)$, where the functional
$\mathcal F_{\varepsilon,H}$ is introduced in (\ref{V-lem-psi<u}).
Let us recall the constant $a_0\in]0,1[$ introduced in
Theorem~\ref{corollary}.

\begin{prop}\label{upperbound}
Let $(\psi,A)$ be a minimizer of (\ref{V-EGL}) and
$\varphi=\displaystyle\frac{\psi}{u_\varepsilon}$. Assume that
$a\in]0,a_0[$. There exist constants $C_*>0$, $\varepsilon_0>0$ such
that, when the applied magnetic field $H$ satisfies (\ref{H=Hc}),
the following upper bound of the energy holds,
$$\mathcal F_{\varepsilon}(\varphi,A)\leq
H^2J_0(\varepsilon)+(C_*-\lambda)(\ln|\ln\varepsilon|)^2,
\quad\forall~\varepsilon\in]0,\varepsilon_0]\,.$$
\end{prop}

The proof of Proposition~\ref{upperbound} is by constructing a
suitable test configuration having vortices and by computing its energy. The estimate
of the energy of the test configuration relies on a careful analysis
of a Green's potential.

\subsection{Analysis of a Green's potential}
This section is devoted to an analysis of  a Green's kernel, i.e. a
fundamental solution of the differential operator\break $-{\rm
  div}\,\left(\displaystyle\frac1{u_\varepsilon^2(x)}\nabla\right)+1$.
The existence and the properties of this function, taken from
\cite{AfSaSe, Stamp}, are given in the next lemma.

\begin{lem}\label{Green}
For every $y\in\Omega$ and $\varepsilon\in]0,1]$, there exists a
symmetric function
$\overline\Omega\times\overline\Omega\ni(x,y)\mapsto
G_\varepsilon(x,y)\in\mathbb R_+$ such that~:
\begin{equation}\label{equation-Green}
\left\{\begin{array}{rl} -{\rm
div}\,\left(\displaystyle\frac1{u_\varepsilon^2(x)}\nabla_x
G_\varepsilon(x,y)\right)+G_\varepsilon(x,y)=\delta_y(x)&{\rm
in}~\mathcal
D'(\Omega),\\
G_\varepsilon(x,y)\big{|}_{x\in\partial\Omega}=0.
\end{array}\right.
\end{equation}
Moreover, $G_\varepsilon$ satisfies the following properties:
\begin{enumerate}
\item There exists a constant $C>0$ such that
$$G_\varepsilon(x,y)\leq C\left(|\,\ln|x-y|\,|+1\right)\,,\quad
\forall~(x,y)\in\overline\Omega\times\overline\Omega\setminus\Delta\,,~
\forall~\varepsilon\in]0,1]\,,$$
where $\Delta$ denotes the diagonal in $\mathbb R^2$.
\item The function $v_\varepsilon:
\overline\Omega\times\overline\Omega\ni (x,y)\mapsto
G_\varepsilon(x,y)+\displaystyle\frac{u_\varepsilon^2(x)}{2\pi}\ln|x-y|$
is in the class $C^1(\overline\Omega\times\overline\Omega\,;\mathbb
R)$\,.
\item Given $q\in[1,2[$, there exists a constant $C>0$ such that
$$\|v_\varepsilon(\cdot,y)\|_{W^{1,q}(\Omega)} \leq C,
\quad \forall~y\in\overline\Omega,\quad\forall~\varepsilon\in]0,1]\,.$$
\item For any compact set $K\subset\Omega$, there exist  constants
$ C>0$ and $\varepsilon_0>0$ such that,
$\forall~\varepsilon\in]0,\varepsilon_0]$,
$$\left|G_\varepsilon(x,y)
+\frac{u_\varepsilon^2(x)}{2\pi}\ln|x-y|\,\right|\leq C\left\|
\nabla u_\varepsilon(x)\right\|_{L^\infty(K)} ,\quad\forall~x\in
K,~\forall~y\in\overline\Omega.$$
\end{enumerate}
\end{lem}

\begin{corol}\label{corl:Holdernorm}
Assume that $a\in]0,a_0[$ and $R'\in]R,1[$.
There exist constants
$C>0$, $\alpha\in]0,1[$ and $\varepsilon_0>0$ such that, for all
$\varepsilon\in]0,\varepsilon_0[$ and $2(R_\varepsilon-R)<
\eta(\varepsilon)<1$, we have
$$\left\|v_\varepsilon(\cdot,y)\right\|
_{C^{0,\alpha}( D(0,R')\setminus D(0,R+\eta(\varepsilon)))}\leq
\frac{C}{\eta(\varepsilon)^2}\,,\quad
\forall~y\in\overline\Omega\,.$$
Here
$$v_\varepsilon(x,y)=G_\varepsilon(x,y)
+\frac{u_\varepsilon^2(x)}{2\pi}\ln|x-y|\,.$$
\end{corol}
\paragraph{\bf Proof.}
Let $\chi\in C_0^\infty(\mathbb R;\mathbb R)$ be a cut-off function such that
$$0\leq\chi\leq1\quad{\rm in~}\mathbb R\,,\quad
\chi\equiv 1\quad {\rm in~}[1,\infty[\,,\quad
\chi\equiv0\quad{\rm in}~]-\infty,\frac12[\,.$$
Set
$$\chi_\eta(x)=\chi\left(\frac{|x|}{\eta}\right)\,,\quad
\widetilde v_\varepsilon(x)=\chi_\eta(x)\,v_\varepsilon(x)\,,\quad
\forall~x\in D(0,1)\,.$$
The function $\widetilde v_\varepsilon$ satisfies the equation
$$-{\rm div}\left(\frac1{u_\varepsilon^2}\nabla_x \widetilde
v_\varepsilon\right)
+\widetilde v_\varepsilon(x,y)
=\chi_\eta\,f_\varepsilon(x,y)+w_\varepsilon(x,y)\,,$$
where
$$f_\varepsilon(x,y)=\frac{u_\varepsilon^2(y)}{\pi u_\varepsilon^3(x)}
\nabla u_\varepsilon(x)\cdot\nabla_x\ln|x-y|
-\frac{u_\varepsilon^2(y)}{2\pi}\ln|x-y|\,,$$
and
$$w_\varepsilon(x,y)
=\frac{\nabla_x
  v_\varepsilon(x,y)\cdot\nabla\chi_\eta(x)}
{u_\varepsilon^2(x)}-2\frac{v_\varepsilon(x,y)}{u_\varepsilon^3(x)}
\nabla u_\varepsilon(x)\cdot\nabla
\chi_\eta(x)+\frac{v_\varepsilon(x,y)}{u_\varepsilon^2(x)}
\Delta\chi_\eta(x)\,.$$
Let us notice also that it results from
Theorems~\ref{V-thm-kach3}~and~\ref{corollary}
$$\|\nabla
u_\varepsilon\|_{L^\infty( D(0,1)\setminus D(0,R+\eta))}\leq C,\quad
 u_\varepsilon\geq \sqrt{a}\quad{\rm in }~ D(0,1)\,.$$
Thanks to the above properties of the function $u_\varepsilon$,
we deduce that for a given $q\in[1,2[$, there exists a constant $C>0$
such that
$$\|\chi_\eta\,f_\varepsilon(\cdot,y)\|_{L^q(D(0,1))}\leq
C\,,\quad\forall
~y\in\overline\Omega\,,~\forall~\varepsilon\in]0,1]\,.$$
On the other hand, for a given $q\in[1,2[$, it is known that the function
$v_\varepsilon(\cdot,y)$ is bounded in $W^{1,q}(\Omega)$ (see
Lemma~\ref{Green}).
Thus, we get the following estimate for the functions
$\widetilde v_\varepsilon$ and $w_\varepsilon$:
$$\|\widetilde v_\varepsilon(\cdot,y)\|_{W^{1,q}(D(0,1))}\leq C\,,\quad
\|w_\varepsilon(\cdot,y)\|_{L^q(D(0,1))}\leq \frac{C}{\eta^2}\,,
~\forall~y\in\overline{D(0,1)}\,,~\forall~\varepsilon\in]0,1]\,.$$
Let $R'\in]R,1[$. Thanks to the equation of $\widetilde v_\varepsilon$,
Theorem~2 of \cite{Me} implies that there exist $p>2$ and $p'<2$ such
that
\begin{equation}\label{Meyers1}
\|\nabla \widetilde v_\varepsilon(\cdot,y)\|_{L^p(D(0,R'))}
\leq C\left(\|\nabla \widetilde v_\varepsilon(\cdot,y)\|_{L^{p'}(D(0,1))}
+\|w_\varepsilon(\cdot,y)+
\chi_\eta f_\varepsilon(\cdot,y)\|_{W^{-1,p}(D(0,1))}\right)\,.\end{equation}
We may choose $q\in]1,2[$ such that $W^{-1,p}\subset L^q$ and $p'<q$. Thus,
the above estimate reads as:
$$\|\nabla \widetilde v_\varepsilon(\cdot,y)\|_{L^p(D(0,R'))}\le
\frac{C}{\eta^2}\,,\quad
 \forall~y\in\overline{D(0,1)}\,.$$
Since $\widetilde v_\varepsilon$ is bounded uniformly in
$W^{1,q}(\Omega)$
(see Lemma~\ref{Green}), we get by  Poincar\'e's inequality:
$$\|\widetilde
v_\varepsilon(\cdot,y)\|_{W^{1,p}(D(0,R'))}\leq\frac{C}{\eta^2}\,,\quad
 \forall~y\in\overline{D(0,1)}\,.$$
Since $p>2$, the Sobolev embedding theorem yields the bound
$$\|\widetilde
v_\varepsilon(\cdot,y)\|_{C^{0,\alpha}(D(0,R'))}\leq\frac{C}{\eta^2}\,,\quad
 \forall~y\in\overline{D(0,1)}\,,$$
for some $\alpha\in]0,1[$. This estimate is nothing but the result of
    the corollary once  we
remember the definition of the function  $\widetilde
v_\varepsilon$\,.
\hfill$\Box$\\

The next  lemma provides us with  points enjoying
useful properties. These points will serve to be the centers of the
vortices of the test configuration that we shall construct in the next section.

\begin{lem}\label{choice-ai}
There exist   constants $K>0$, $c\in]0,1[$, and for each
$\varepsilon\in]0,1[$ and $n(\varepsilon)\in\mathbb
N\,\cap\,[1,\frac{c}2\,\varepsilon^{-1}[$\,, there exist points
$(a_i)_{i=1}^{n(\varepsilon)}\subset
\partial D(0,r_\varepsilon)$ and $\delta(\varepsilon)\in]0,1[$ such that
\begin{eqnarray*}
&&\delta(\varepsilon)\ll1\quad{\rm as~}\varepsilon\to0\,\\
&&\frac{c}{ n(\varepsilon)}\leq |a_{i+1}-a_i|\leq \delta(\varepsilon)+\frac{c}{n(\varepsilon)}\,,\quad |v_\varepsilon(a_i,a_i)|\leq K\ln|\ln\varepsilon|\,,\\
&&\forall~i\in\{1,2,\cdots,n(\varepsilon)\}\,,\quad\forall~\varepsilon
\in]0,1]\,.
\end{eqnarray*}
Here the function $v_\varepsilon$ has been introduced in
Lemma~\ref{Green}\,,
and $r_\varepsilon=R
+\displaystyle\frac{\ln|\ln\varepsilon|}{|\ln\varepsilon|}$\,.
\end{lem}
\paragraph{\bf Proof.}
The proof is actually due to the following bound
\begin{equation}\label{mainbound}
\int_{\partial D(0,r_\varepsilon)}|v_\varepsilon(x,x)|\,\md x\leq
C\ln|\ln\varepsilon|\,,\end{equation} that holds uniformly in $\varepsilon\in]0,1]$.
Let us show why this bound holds. We cover $\partial
D(0,r_\varepsilon)$ by $\mathcal N$ balls $(B(y_i,\zeta))_i$, with
$(y_i)_i\subset
\partial D(0,r_\varepsilon)$ and $\zeta\in]0,1[$ is to be chosen
appropriately. We introduce a scaled partition of unity
$\chi_i^\zeta$ such that
$$\sum_{i=1}^{\mathcal N}|\chi_i^\zeta|=1~{\rm in~}\partial
D(0,r_\varepsilon)\,,\quad {\rm supp}\,\chi_i^\zeta\subset
B(y_i,\zeta)\,,\quad\forall~i\in\{1,\cdots,\mathcal N\}\,.$$ Then
\begin{equation}\label{IMS}
\int_{\partial D(0,r_\varepsilon)} |v_\varepsilon(x,x)|\,\md x=
\sum_{i=1}^{\mathcal N}\int_{\partial D(0,r_\varepsilon)}
|\chi_i^\zeta(x)\,v_\varepsilon(x,x)|\,\md x\,.
\end{equation}
By Corollary~\ref{corl:Holdernorm}, we write for all
$i\in\{1,\cdots,\mathcal N\}$\,
\begin{eqnarray*}
\int_{\partial D(0,r_\varepsilon)}
|\chi_i^\zeta(x)\,v_\varepsilon(x,x)|\,\md x&\leq& \int_{\partial
D(0,r_\varepsilon)}|\chi_i^\zeta(x)\,
v_\varepsilon(x,y_i)|\,\md\mu_*(x)\\
&& +\frac{C}{\eta^2}\int_{\partial
D(0,r_\varepsilon)}|\chi_i^\zeta(x)|\,|x-y_i|^\alpha\,
\md\mu_*(x)\,,\end{eqnarray*} where $\alpha\in]0,1[$,
$\eta=R-r_\varepsilon$, and $\mu_*$ is the Lebesgue measure in
$\partial D(0,r_\varepsilon)$.\\
Recalling that ${\rm supp}\,\chi_i^\zeta\subset B(y_i,\zeta)$, we
get upon choosing $\zeta=\eta^{2/\alpha}$ and summing up over $i$,
$$\sum_{i=1}^{\mathcal N}\frac{1}{\eta^2}\int_{\partial
D(0,r_\varepsilon)}|\chi_i^\zeta(x)|\,|x-y_i|^\alpha\,
\md\mu_*(x)\leq C\,.$$ 
On the other hand, by Lemma~\ref{Green}, there exists a constant $C>0$
such that
\begin{eqnarray*} 
\int_{\partial D(0,r_\varepsilon)}|\chi_i^\zeta(x)\,
v_\varepsilon(x,y_i)|\,\md\mu_*(x)&\leq&
C\int_{B(y_i,\zeta)\cap \partial D(0,r_\varepsilon)}|\,\ln|x-y_i|\,
|\,\md\mu_*(x)\\
&\leq&C\zeta |\ln\zeta|\,.\end{eqnarray*}
Recalling our choice of  $\zeta=\eta^{2/\alpha}$ and $\eta=
\mathcal O(|\ln\varepsilon|^{-1/2})$, and summing up over $i$,
we get
for a new constant $C>0$
$$\sum_{i=1}^{\mathcal N}\int_{\partial D(0,r_\varepsilon)}
|\chi_i^\zeta(x)\,v_\varepsilon(x,x)|\,\md x\leq C\mathcal
N\times\zeta\,
\ln|\ln\varepsilon|\leq C\ln|\ln\varepsilon|\,,$$
where we have used that $\mathcal N\times 2\pi \zeta\approx 2\pi
r_\varepsilon\to 2\pi R$. Substituting in (\ref{IMS}), we obtain
the desired bound  (\ref{mainbound}).\\
Now, defining the function
$$f_\varepsilon(t):[0,1[\ni t\mapsto
|v_\varepsilon\left( r_\varepsilon e^{2\pi i \,t},r_\varepsilon e^{2\pi
  i\,t}
\right)|\,,$$
and applying Lemma~\ref{integration} below, we get the desired sequence of
points.\hfill$\Box$\\

\begin{lem}\label{integration}
Let $(f_\varepsilon)_{\varepsilon\in]0,1]}
\subset C([0,1],\mathbb R_+)$ be a family of
continuous functions. Assume that  there exists a constant $C>0$
such that
$$\|f_\varepsilon\|_{L^1([0,1])}\leq C\ln|\ln\varepsilon|\,,\quad\forall~\varepsilon\in
]0,1]\,.$$
There exist
constants $K>0$ and $c_0\in]0,1[$ such that, given a
family $(N(\varepsilon))\subset\mathbb N$ satisfying
$N(\varepsilon)\gg1$,  there exists
a family $(\delta(\varepsilon))\subset ]0,1[$ and a sequence
$(t_{m}^\varepsilon)_{m\in\mathbb N}\subset ]0,1[$
and
\begin{eqnarray*}
&&|f_\varepsilon(t_i)|\leq K\ln|\ln\varepsilon|\,,
\quad \frac{c_0}{N(\varepsilon)}
\leq\left|t_{i+1}^\varepsilon-t_i^\varepsilon\right|\leq
\delta(\varepsilon)+\frac{c_0}{N(\varepsilon)}\,,\\
&&\forall~i\in\{1,2,\cdots,N(\varepsilon)\}\,,\quad\forall~\varepsilon\in
]0,1]\,.
\end{eqnarray*}
\end{lem}
\paragraph{\bf Proof.}
Let us introduce, for a given $K>0$, the set
$$E_K^\varepsilon=\{t\in[0,1]~:~|f_\varepsilon(t)|<K\ln|\ln\varepsilon|\}\,.$$
Using the uniform bound on  $\|f_\varepsilon\|_{L^1([0,1])}$,
we can choose
  $K$ sufficiently large such that
$$|E_K^\varepsilon|\geq\frac12\quad\forall~\varepsilon\in]0,1]\,,$$
where $|\cdot|$ denotes the Lebesgue measure.\\
Let $\varepsilon\in]0,1]$.
Since the function $f_\varepsilon$ is continuous,
the set $E_K^\varepsilon$ is open.
Thus, we essentially meet two cases: Either
there exists an interval
$$]x_\varepsilon-\delta_\varepsilon,x_\varepsilon+\delta_\varepsilon[\subset
    E_K^\varepsilon$$
with $\displaystyle
\liminf_{\varepsilon\to0}\delta_\varepsilon>0$ (in which case the
    statement of the theorem becomes evidently true), or there exists
a constant $c_0\in]0,\frac12[$  and  possibly infinitely many disjoint
intervals
$$\bigcup_{i\in I_\varepsilon}]x_i^\varepsilon-\delta_i^\varepsilon,
x_i^\varepsilon+\delta_i^\varepsilon[\subset E_K^\varepsilon$$
such that $x_1^\varepsilon<x_2^\varepsilon<\cdots$ and
\begin{equation}\label{shape-E}
\limsup_{\varepsilon\to0}\left(\sum_{i\in I_\varepsilon}
\Delta_i
\right)=0\,,\quad
\liminf_{\varepsilon\to0}\left(\sum_{i\in I_\varepsilon}
\delta_i^\varepsilon
\right)\geq \frac{c_0}2\,,\end{equation}
where $\Delta_i=\left|x_{i+1}^\varepsilon-\delta_{i+1}^\varepsilon
-x_i^\varepsilon-\delta_i^\varepsilon\right|$.\\
Consequently, setting    $t_1^\varepsilon=x_1^\varepsilon$
 and $\mathcal E_K^\varepsilon=\displaystyle
\bigcup_{i\in I_\varepsilon}]x_i^\varepsilon-\delta_i^\varepsilon,
x_i^\varepsilon+\delta_i^\varepsilon[$, we get
$$\mathcal E_K^\varepsilon{\bigg \backslash}
\left]t_1^\varepsilon,t_1^\varepsilon+\frac{c_0}{8N(\varepsilon)}\right[
\not=\emptyset\,.$$
So, setting
$$t_2^\varepsilon=\inf\left(\mathcal E_K^\varepsilon{\bigg \backslash}
\left]t_1^\varepsilon,t_1^\varepsilon+\frac{c_0}{8N(\varepsilon)}\right[\,
\right)
\,,$$
we get, thanks in particular to (\ref{shape-E}),
$$
t_2^\varepsilon-t_1^\varepsilon\leq 2\sum_i\Delta_i+
\frac{c_0}{8N(\varepsilon)}\,,$$
and
$$\mathcal E_K^\varepsilon{\bigg \backslash}
\left]t_1^\varepsilon,t_2^\varepsilon+\frac{c_0}{8N(\varepsilon)}\right[
\not=\emptyset\,.$$
Therefore, we set
$$t_3^\varepsilon= \inf\left(\mathcal E_K^\varepsilon{\bigg \backslash}
\left]t_1^\varepsilon,
t_2^\varepsilon+\frac{c_0}{8N(\varepsilon)}\right[\,\right)
\,.$$
By induction, given $n\leq N(\varepsilon)$, we can construct points
$t_2^\varepsilon<t_3^\varepsilon<\cdots<t_n^\varepsilon$  such that
$$
\frac{c_0}{8N(\varepsilon)}\leq
t_{i+1}^\varepsilon-t_i^\varepsilon\leq\frac{c_0}{8N(\varepsilon)}
+\sum_i\Delta_i\,,\quad\forall~i\in\{1,\cdots,n\}\,,$$
yielding therefore the desired sequence $(t_m^\varepsilon)$ with
$\delta(\varepsilon)=\displaystyle
\sum_{i\in I_\varepsilon}\Delta_i$.\hfill$\Box$\\

\subsection{The test configuration}
We know from Theorem~\ref{corollary} that the function
$\xi_\varepsilon(x)$ achieves its unique minimum on the circle
$\partial D(0,R_\varepsilon)$ with $\varepsilon\ll
R_\varepsilon-R\ll\varepsilon^\alpha$, for $\alpha\in]0,1[$.\\
Since we  expect vortices of a minimizer of (\ref{V-EGL}) to be
pinned on the circle $\partial D(0,R_\varepsilon)$, and to be
uniformly distributed along this circle, we  construct a test
configuration
whose vortices are placed, for technical reasons, on the circle
$\partial D(0,r_\varepsilon)$, with $r_\varepsilon=R+\displaystyle
\frac{\ln|\ln\varepsilon|}{|\ln\varepsilon|}$. We mention that
similar constructions have been also introduced
in the papers \cite{AfAlBr, AlBr06, AlBr}.\\
Let $n(\varepsilon)\in\mathbb
N\,\cap\,]1,\frac{c}2\varepsilon^{-1}[$
for an appropriate constant $c\in]0,1[$.
Lemma~\ref{choice-ai} provides us with $n(\varepsilon)$  points
$(a_i)_{i=1}^{n(\varepsilon)}$ on the circle $\partial
D(0,r_\varepsilon)$,  that satisfy in particular
$$\overline{B(a_i,\varepsilon)}\cap \overline{B(a_j,\varepsilon)}
=\emptyset\,,\quad
\forall~i\not =j\,.$$
We define a measure $\mu$ by:
\begin{equation}\label{measure}
\mu(x)=\left\{\begin{array}{lll}
0&{\rm if}&x\not\in\cup_{i}B(a_i,\varepsilon)\\
\displaystyle\frac2{\varepsilon^2}&{\rm if}&x\in\overline{\cup_i
B(a_i,\varepsilon)},
\end{array}\right.\end{equation}
and a function $h'$ in $\Omega=D(0,1)$ by
\begin{equation}\label{V-Def-f'}
\left\{
\begin{array}{rl}
-{\rm div}\,\left(\displaystyle\frac1{u_\varepsilon^2}\nabla
h'\right)
+h'=\mu&{\rm in}~\Omega,\\
h'=0&{\rm on}~\partial \Omega.\end{array}\right.\end{equation} We
notice that
$$\int_{B(a_i,\varepsilon)}\mu\,\md x=2\pi,\quad\forall~i=1,2,\cdots,
n(\varepsilon)\,,\quad
\int_{\mathbb R^2}\mu\,\md x=2\pi\,n(\varepsilon)_,.$$
We define  an induced magnetic field $h=h'+h_\varepsilon$  (here
$h_\varepsilon$ has been introduced in (\ref{V-hepsilon'})). Then we
define an induced magnetic potential
$A=A'+\frac{H}{u_\varepsilon^2}\nabla^\bot h_\varepsilon$ by
  taking simply
$${\rm curl}\,A'=h'.$$
This choice is always possible as one can take $A'=\nabla^\bot g$
with $g\in H^2(\Omega)$ such that
$\Delta g=h'$.\\
We turn now to define an order parameter $\psi$ which we take in the
form
\begin{equation}\label{V-vsol-psi'}
\psi=u\,u_\varepsilon=\rho\,e^{i\phi}\,u_\varepsilon,\end{equation}
where $\rho$ is defined by:
\begin{equation}\label{V-vsol-rho'}
\rho(x)=\left\{
\begin{array}{cll}
0&{\rm if}&x\in \cup_iB(a_i,\varepsilon),\\
1&{\rm if}&x\not\in\cup_iB(a_i,2\varepsilon),\\
\displaystyle\frac{|x-a_i|}\varepsilon-1&{\rm if}&\exists\,i~{\rm
s.t.}~x\in B(a_i,2\varepsilon) \setminus B(a_i,\varepsilon).
\end{array}\right.
\end{equation}
The phase $\phi$ is defined (modulo $2\pi$) by the relation:
\begin{equation}\label{V-vsol-phi'}
\nabla\phi-A'=-\frac{1}{u_\varepsilon^2}\nabla^\bot h'\quad{\rm in}~
\Omega\setminus \cup_iB(a_i,\varepsilon),
\end{equation}
and we emphasize here that we do not need to define $\phi$ in
regions where $\rho$ vanishes.\\

\begin{lem}\label{int-Green}
There exist constants $\varepsilon_0\in]0,1[$ and $C>0$ such that
\begin{eqnarray*}
&&\hskip-0.5cm\int_\Omega\left(\frac{1}{u_\varepsilon^2(x)}|\nabla
h'|^2+|h'|^2\right)\,\md
x\,\md y\\&&\leq
2\pi
u_\varepsilon^2(r_\varepsilon)\,n(\varepsilon)|\ln\varepsilon|+C\,n(\varepsilon)\ln|\ln\varepsilon|+C\,
[n(\varepsilon)]^2+o\left([n(\varepsilon)]^2\right).
\end{eqnarray*}
\end{lem}
\paragraph{\bf Proof.}
Notice that the field $h'$ can be expressed by means of the function
$G_\varepsilon$ introduced in Lemma~\ref{Green},
\begin{equation}\label{h'-Green}
h'(x)=\int_\Omega G_\varepsilon(x,y)\,\mu(y)\,\md
y\,,\quad\forall~x\in\Omega.
\end{equation}
Therefore, we get the identity
\begin{equation}\label{Appendix-lem}
\int_\Omega\left(\frac{1}{u_\varepsilon^2(x)}|\nabla
h'|^2+|h'|^2\right)\,\md x=
\int_{\Omega\times\Omega}G_\varepsilon(x,y)\,\mu(x)\mu(y)\,\md
x\,\md y\,,\end{equation}
which shows that it is sufficient  to estimate $\displaystyle\int_{\Omega\times\Omega}
G_\varepsilon(x,y)\,\mu(x)\mu(y)\,\md x\,\md y$.
We decompose the integral  $\displaystyle\int_{\Omega\times\Omega}
G_\varepsilon(x,y)\,\mu(x)\mu(y)\,\md x\,\md y$ into two terms:
\begin{eqnarray}\label{decomposition:integral}
&&\int_{\Omega\times\Omega} G_\varepsilon(x,y)\,\mu(x)\mu(y)\,\md
x\,\md y=\\
&&\hskip1cm\sum_{i\not=j} \int_{B(a_i,\varepsilon)\times
B(a_j,\varepsilon)} G_\varepsilon(x,y)\,\mu(x)\mu(y)\,\md
x\,\md y\nonumber\\
&&\hskip1cm+\sum_{i=1}^{n(\varepsilon)}
\int_{B(a_i,\varepsilon)\times B(a_i,\varepsilon)}
G_\varepsilon(x,y)\,\mu(x)\mu(y)\,\md x\,\md y\,.\nonumber
\end{eqnarray}
Let us estimate the first term.  We write using
Lemma~\ref{Green},
\begin{eqnarray}\label{ub-proof-1}
&&\hskip-0.5cm\sum_{i\not=j} \int_{B(a_i,\varepsilon)\times
B(a_j,\varepsilon)} G_\varepsilon(x,y)\,\mu(x)\mu(y)\,\md x\,\md
y\nonumber\\
&&\leq C\sum_{i\not=j} \int_{B(a_i,\varepsilon)\times
B(a_j,\varepsilon)}
\big{(}\left|\,\ln|x-y|\,\right|+1\big{)}\mu(x)\mu(y)\,\md x\md y\,.
\end{eqnarray}
Now, recalling the definition of $\mu$ in (\ref{measure}) and the
choice of the points $(a_i)$ in Lemma~\ref{choice-ai} , we get
\begin{equation}\label{ub-proof1*}
\sum_{i\not=j}\frac{c}{n(\varepsilon)^2}\int_{B(a_i,\varepsilon)\times
B(a_j,\varepsilon)}
|\,\ln|x-y|\,|\,\mu(x)\mu(y)\md x\md y\leq
C\,,
\end{equation}
where $C>0$ is any constant such that
$$C>\displaystyle\int_{\partial D(0,R)\times \partial
D(0,R)}\left|\,\ln|x-y|\,\right|\,
\md \mu_*(x)\md \mu_*(y)$$ and $\md \mu_*$
is the arc-length measure on the circle $\partial D(0,R)$.\\
Therefore, (\ref{ub-proof-1}) becomes for a
new constant $C>0$,
\begin{equation}\label{ub-proof1-conc}
\sum_{i\not=j}
\int_{B(a_i,\varepsilon)\times
B(a_j,\varepsilon)}G_\varepsilon(x,y)\,\mu(x)\mu(y)\,\md x\md y\leq
C\, n(\varepsilon)^2\,.
\end{equation}
Again, using Corollary~\ref{corl:Holdernorm}, we estimate
\begin{eqnarray}\label{ub-proof2:decomp}
&&\hskip-0.5cm\int_{B(a_i,\varepsilon)\times B(a_i,\varepsilon)}
G_\varepsilon(x,y)\,\mu(x)\mu(y)\,\md
x\,\md y\\
&&=\frac{4}{\varepsilon^4}
\int_{B(a_i,\varepsilon)\times
B(a_i,\varepsilon)}
\left(\frac{u_\varepsilon^2(x)}{2\pi}\ln\frac1{|x-y|}
+|v_\varepsilon(x,y)|\right)\,\md x\md y\,.\nonumber
\end{eqnarray}
On the one hand, we have
\begin{eqnarray*}
&&\hskip-0.5cm
\frac{4}{\varepsilon^4}
\int_{B(a_i,\varepsilon)\times
B(a_i,\varepsilon)}
\frac{u_\varepsilon^2(x)}{2\pi}\ln\frac1{|x-y|}\,\md x\md y\\
&&= 4\int_{B(0,1)\times B(0,1)}
\frac{u_\varepsilon^2(a_i+\frac{z'}\varepsilon)}{2\pi}
\ln\left[\frac1{\varepsilon |z'-w'|}\right]\,\md
z'\md w'\,.
\end{eqnarray*}
Recall that the function $u_\varepsilon$ is radial and that
$|a_i|=r_\varepsilon
$. By Theorem~\ref{corollary}, we have
$$
\left|u_\varepsilon^2\left(a_i+\frac{z'}\varepsilon\right)
-u_\varepsilon^2(r_\varepsilon)\right|\leq \mathcal
O(\varepsilon),\quad\forall~z'\in B(0,1)\,.$$
Therefore,
\begin{eqnarray}\label{ub-proof-3*}
&&\hskip-0.5cm
\int_{B(a_i,\varepsilon)\times B(a_i,\varepsilon)}
\frac{u_\varepsilon^2(x)}{2\pi}\ln\frac1{|x-y|}\,\mu(x)\mu(y)\,\md
x\,\md y\nonumber\\
&&\leq 4\int_{B(0,1)\times B(0,1)}
\frac{u_\varepsilon^2(r_\varepsilon)+\mathcal O(\varepsilon)}{2\pi}
\ln\left[\frac1{\varepsilon |z'-w'|}\right]
\,\md z'\md w'\nonumber\\
&&=2\pi\,u_\varepsilon^2(r_\varepsilon)|\ln\varepsilon|+o(1)\,.
\end{eqnarray}
On the other hand, assuming that the following estimate holds
\begin{equation}\label{ub-proof2-claim}
\limsup_{\varepsilon\to0}\frac1{n(\varepsilon)\ln|\ln\varepsilon|}\sum_{i=1}
^{n(\varepsilon)}
\frac{4}{\varepsilon^4}\int_{B(a_i,\varepsilon)\times
  B(a_i,\varepsilon)}|v_\varepsilon(x,y)|\,\md
x\md y
\leq C\,,
\end{equation}
then (\ref{ub-proof2:decomp}) becomes
\begin{equation}\label{ub-proof-3}
\sum_{i=1}^{n(\varepsilon)}
\int_{B(a_i,\varepsilon)\times
  B(a_i,\varepsilon)}G_\varepsilon(x,y)\,\mu(x)\mu(y)\,\md x\md y\leq
\left[u_\varepsilon^2(r_\varepsilon)|\ln\varepsilon|+C\ln|\ln\varepsilon|
\right]
n(\varepsilon)\,.
\end{equation}
Combining (\ref{ub-proof-1}) and
(\ref{ub-proof-3}), and using (\ref{Appendix-lem}), we get the result of Lemma~\ref{int-Green}.\\
It remains to prove the claim in
(\ref{ub-proof2-claim}).\\
{\it Proof of (\ref{ub-proof2-claim}).}\\
By Corollary~\ref{corl:Holdernorm}, we write,
$$\frac{4}{\varepsilon^4}
\int_{B(a_i,\varepsilon)\times B(a_i,\varepsilon)}
|v_\varepsilon(x,y)|\,\md x\md y\leq
4\pi|v_\varepsilon(a_i,a_i)|+\frac{C\varepsilon^\alpha}{\eta^2}\,,$$
where $\alpha\in]0,1[$ and $\eta=r_\varepsilon-R$.
Using our particular choice of
$\eta=\mathcal O(|\ln\varepsilon|^{-1/2})$, we deduce that
$$\frac1{n(\varepsilon)}\sum_{i=1}^{n(\varepsilon)}
\frac{4}{\varepsilon^4}
\int_{B(a_i,\varepsilon)\times B(a_i,\varepsilon)}
|v_\varepsilon(x,y)|\,\md x\md y\leq  \frac{4\pi}{n(\varepsilon)}
\left(\sum_{i=1}^{n(\varepsilon)}|v_\varepsilon(a_i,a_i)|+o(1)\right)\,.$$
Recalling the choice of the points $(a_i)$ in  Lemma~\ref{choice-ai},
we  see that the right hand side above
is uniformly bounded by a constant times $\ln|\ln\varepsilon|$, 
yielding the result in (\ref{ub-proof2-claim}).
\hfill$\Box$\\

In the next lemma, we state a decomposition of the energy due to
\cite{BeRi}.

\begin{lem}\label{Beth-Riv}
Consider $(u,A)\in H^1(\Omega;\mathbb C)\times H^1(\Omega;\mathbb
R^2)$ and define
$$A'=A-\frac{H}{u_\varepsilon^2}\nabla^\bot h_\varepsilon,$$
where $u_\varepsilon$ and $h_\varepsilon$ are introduced in
Theorem~\ref{V-thm-kach3} and (\ref{V-hepsilon'}) respectively. Then
we have the decomposition of the energy,
\begin{align*}
\begin{split}
\mathcal F_{\varepsilon,H}(u,A)=&H^2J_0(\varepsilon)
+\int_\Omega\left( u_\varepsilon^2|(\nabla-iA')u|^2+|{\rm
curl}\,A'|^2+\frac1{\varepsilon^2}
u_\varepsilon^4(1-|u|^2)^2\right)\md x
\\&+
2H\int_\Omega (h_\varepsilon-1)\bigg{[}{\rm curl}\big{(}
A'+(iu,\nabla_{A'}u)\big{)}\bigg{]}\md x\\
&+H^2\int_\Omega\frac1{u_\varepsilon^2}\left(|u|^2-1\right)
|\nabla h_\varepsilon|^2\,\md x.
\end{split}
\end{align*}
Here, the functional $\mathcal F_{\varepsilon,H}$ and the energy
$J_0(\varepsilon)$ are introduced in (\ref{V-reducedfunctional}) and
(\ref{V-J0}) respectively. \end{lem}

\paragraph{\bf Proof of Proposition~\ref{upperbound}.}
Let $(\psi,A)$ be the test configuration constructed in
(\ref{V-Def-f'})-(\ref{V-vsol-phi'}), and put
$\varphi=\frac\psi{u_\varepsilon}$. By Lemma~\ref{V-lem-psi<u}, it is
sufficient to establish the upper bound
$$\mathcal F_{\varepsilon,H}(\varphi,A)\leq
H^2J_0(\varepsilon)+(C_*-\lambda)(\ln|\ln\varepsilon|)^2\,.
$$
By the construction of $\varphi$ and Theorem~\ref{V-thm-kach3}, we
get,
$$\frac1{\varepsilon^2}
\int_\Omega u_\varepsilon^2(1-|\varphi|)^2\,\md x=\mathcal O(1)\,.$$
By Lemma~\ref{V-lem-controlgrad}, we have
$$\int_\Omega\frac1{u_\varepsilon^2}\left(|\varphi|^2-1\right)
|\nabla h_\varepsilon|^2\,\md x\leq C\varepsilon^2\,.$$
Let $\mu(\varphi,A')={\rm curl}\big{(}
A'+(i\varphi,\nabla_{A'}\varphi)\big{)}$. Notice that
$$\mu(\varphi,A')=\left\{\begin{array}{l}
0\quad {\rm in}~\Omega\setminus\bigcup_i B(a_i,2\varepsilon)\,,\\
\mu+\mu_r(\varphi,A')\quad
{\rm in}~\bigcup_iB(a_i,2\varepsilon)\,,
\end{array}\right.$$
where $\mu$ is the measure defined in (\ref{measure}) and
$$\mu_r(\varphi,A')=
-(|\varphi|^2-1){\rm div}\left(\displaystyle\frac1{u_\varepsilon^2}\nabla
h'\right)+(\nabla^\bot|\varphi|^2)\cdot(\nabla\phi-A).$$
Using the definition of $\mu$ and
Lemma~\ref{V-lem-controlgrad}, we write,
\begin{eqnarray}\label{step1-ub}
&&\hskip-0.5cm2H\int_\Omega (h_\varepsilon-1)\bigg{[}{\rm curl}\big{(}
A'+(i\varphi,\nabla_{A'}\varphi)\big{)}\bigg{]}\md x\nonumber\\
&&=2H(h_\varepsilon(r_\varepsilon)-1)\int_\Omega
\mu(\varphi,A')\,\md x+2H\int_\Omega\left[h_\varepsilon(x)-h_\varepsilon(r_\varepsilon)\right]
\mu(\varphi,A')\,\md x\nonumber\\
&&\leq4\pi n(\varepsilon)(h_\varepsilon(r_\varepsilon)-1)H\\
&&\hskip0.5cm+2H(h_\varepsilon(r_\varepsilon)-1)\int_\Omega\mu_r(\varphi,A')\,\md x+
C\varepsilon\int_\Omega|\mu(\varphi,A')|\,\md x.\nonumber
\end{eqnarray}
Since $|\varphi|=1$ on $\partial B(a_i,2\varepsilon)$, an integration
 by parts yields
\begin{eqnarray*}
&&\hskip-0.5cm\int_{B(a_i,2\varepsilon)}\mu_r(\varphi,A')\,\md x\\
&&=
\int_{B(a_i,2\varepsilon)}\left[
-(|\varphi|^2-1){\rm div}\left(\displaystyle\frac1{u_\varepsilon^2}\nabla
h'\right)+(\nabla^\bot(|\varphi|^2-1))\cdot\frac{\nabla^\bot h'}{u_\varepsilon^2}
  h')\right]\,\md x\\
&&=\int_{B(a_i,2\varepsilon)}\left[
-(|\varphi|^2-1){\rm div}\left(\displaystyle\frac1{u_\varepsilon^2}\nabla
h'\right)+(|\varphi|^2-1){\rm div}\left(\frac1{u_\varepsilon^2}\nabla h'\right)
  h')\right]\,\md x\\
&&=0\,.
\end{eqnarray*}
On the other hand, using the definitions of $\mu$ and $\varphi$,
and Cauchy-Schwarz
inequality, we write,
\begin{eqnarray*}
\int_{B(a_i,2\varepsilon)}|\mu(\varphi,A')|\,\md x&\leq&
\int_{B(a_i,2\varepsilon)}\left(\mu+2|h'|+\frac{C}{\varepsilon}|\nabla
h'|\right)\,\md x\\
&\leq&2\pi +C\left(
\int_{B(a_i,2\varepsilon)}\left(|\nabla h'|^2+|h'|^2\right)\,\md x\right)^{1/2}\,.
\end{eqnarray*}
Therefore, (\ref{step1-ub}) becomes, for a new constant $C>0$,
\begin{eqnarray*}
&&\hskip-0.5cm2
H\int_\Omega (h_\varepsilon-1)\bigg{[}{\rm curl}\big{(}
A'+(i\varphi,\nabla_{A'}\varphi)\big{)}\bigg{]}\md x\\
&&\leq4\pi n(\varepsilon)(h_\varepsilon(r_\varepsilon)-1)H
+C\varepsilon\left[\left(
\int_{\Omega}
\left(|\nabla h'|^2+|h'|^2\right)\,\md x\right)^{1/2}+n(\varepsilon)H\right]\,.
\end{eqnarray*}
Thanks to  Lemma~\ref{Beth-Riv}, we get
\begin{eqnarray*}
\mathcal F_{\varepsilon,H}(\varphi,A)&\leq&
H^2J_0(\varepsilon)+(1+C\varepsilon)
\int_\Omega\left(\frac1{u_\varepsilon^2}|\nabla
h'|^2+|h'-1|^2\right)\md x\\
&&+4\pi n(\varepsilon)(h_\varepsilon(r_\varepsilon)-1)H+\mathcal O
\left(\varepsilon\,n(\varepsilon)H\right)\,.
\end{eqnarray*}
We recall that the magnetic field satisfies
$H=k_\varepsilon|\ln\varepsilon|+\lambda\ln|\ln\varepsilon|$,
and we apply  Lemma~\ref{int-Green} to
deduce the upper bound,
\begin{eqnarray}\label{ub-proof-4}
\mathcal F_{\varepsilon,H}(\varphi,A)&\leq& H^2J_0(\varepsilon)
+2\pi n(\varepsilon)(1-h_\varepsilon(r_\varepsilon))
\left[\left(\frac{u_\varepsilon^2(r_\varepsilon)}{1-
h_\varepsilon(r_\varepsilon)}-2k_\varepsilon\right)|\ln\varepsilon|\right.\nonumber\\
&&\left. +C\frac{n(\varepsilon)+\ln|\ln\varepsilon|}{1-h_\varepsilon(r_\varepsilon)}
-\lambda\ln|\ln\varepsilon|\right]+o(n(\varepsilon)^2+n(\varepsilon)
\sqrt{\varepsilon}\,)\,.
\end{eqnarray}
Recall the definition of $r_\varepsilon=R+\displaystyle
\frac{\ln|\ln\varepsilon|}{|\ln\varepsilon|}$\,.
Thanks to Lemma~\ref{V-lem-controlgrad} and Theorem~\ref{corollary},
we get:
$$\left|\frac{u_\varepsilon^2(r_\varepsilon)}{1-
h_\varepsilon(r_\varepsilon)}-2k_\varepsilon\right|\leq \widetilde C\,
\frac{\ln|\ln\varepsilon|}{|\ln\varepsilon|}\,.$$ Thus, when choosing $C_*>2C+\widetilde C$ and
$n(\varepsilon)=[\,\ln|\ln\varepsilon|\,]$ ($[\,\cdot\,]$ denotes
the largest integer less than $\cdot$), the upper bound
(\ref{ub-proof-4}) becomes,
$$\mathcal F_{\varepsilon,H}(\varphi,A)\leq H^2J_0(\varepsilon)
+(C_*-\lambda)(\ln|\ln\varepsilon|)^2\,,$$ thus achieving the proof of
Proposition~\ref{upperbound}.\hfill$\Box$\\

\section{Proof of main theorems}\label{section:proofs}
\subsection{Proof of Theorem~\ref{thm1}}
Theorem~\ref{LowerBound-thm} provides us with a family of vortex
balls $(B(a_i,r_i))_i$. In particular, when the lower bound of
Corollary~\ref{LowerBound-corol} is  matched with the upper bound
(\ref{V-upperbound}),
permits us to deduce,
$$
0\geq2\pi\,a
\sum_{d_i>0}\left(|\ln\varepsilon|-C\ln|\ln\varepsilon|-k_\varepsilon^{-1}
H\right)|d_i|
+2\pi\sum_{d_i\leq0}\left(a|\ln\varepsilon|-C\ln|\ln\varepsilon|\right)|d_i|\,.
$$
Taking $\lambda_*<-C$ and $\lambda\leq\lambda_*$, we deduce that
$\sum_i|d_i|=0$ whenever the magnetic field satisfies $H<
k_\varepsilon|\ln\varepsilon|+\lambda\ln|\ln\varepsilon|$. 
The energy decomposition of Lemma~\ref{Beth-Riv},
together with Point (4) of Proposition~\ref{V-lem-vortexballs},
yield now the estimate
$$\frac1{\varepsilon^2}\int_\Omega(1-|\varphi|^2)\,\md x\ll1\,$$
which when combined with Lemma~\ref{BBH-thm3.3} gives the desired
result, $|\varphi|\geq\frac12$
in $\overline\Omega$.\\
Now, when $\lambda\geq\lambda_*$, the properties (a)-(c) of
Theorem~\ref{thm1} are consequences of Theorem~\ref{LowerBound-thm}
and Lemma~\ref{upperbound-degree}, which give in particular the upper bound on the total degree
$\sum_i|d_i|\leq C\ln|\ln\varepsilon|$.\\
Assume now that
$H=k_\varepsilon|\ln\varepsilon|+\lambda\ln|\ln\varepsilon|$, with
$\lambda>0$.
When the
lower bound of Corollary~\ref{LowerBound-corol}
\begin{eqnarray*}
\mathcal
F_{\varepsilon,H}(\varphi,A)&\geq& H^2J_0(\varepsilon)+
2\pi a\sum_{d_i>0}\left(|\ln\varepsilon|-C\ln|\ln\varepsilon|-k_\varepsilon^{-1}H\right)d_i\\
&&+2\pi\sum_{d_i\leq0}\left(a|\ln\varepsilon|-C\ln|\ln\varepsilon|\right)|d_i|
\end{eqnarray*}
is matched with the upper bound of Proposition~\ref{upperbound}, we
deduce that
$$2\pi a\sum_{d_i>0}(C-\lambda)\ln|\ln\varepsilon|)d_i\leq
(C_*-\lambda)(\ln|\ln\varepsilon|)^2.
$$
Taking $\mu>\max(C_*,C)$, we deduce the desired lower bound on the
total degree
$$\sum_i|d_i|\geq \sum_{d_i>0}d_i\geq c\ln|\ln\varepsilon|\,.$$
This achieves the proof of Theorem~\ref{thm1}.

\subsection{Proof of Theorem~\ref{thm2}}
Let $(\psi,A)$ be a minimizer of (\ref{V-EGL}) such that $|\psi|>0$ in
$\overline\Omega$. Then
all the degrees
$(d_i)$ in Theorem~\ref{LowerBound-thm} are null:
$$d_i=0\quad\forall~i\,.$$
It results now  from the upper bound $\mathcal
F_{\varepsilon,H}(\varphi,A)\leq H^2 J_0(\varepsilon)$, the lower
bound of Theorem~\ref{LowerBound-thm}
and the energy decomposition of
Lemma~\ref{Beth-Riv}:
$$\int_\Omega\left(|(\nabla-iA')\varphi|^2
+\frac1{2\varepsilon^2}(1-|\varphi|^2)^2+|{\rm
  curl}\,A-Hh_\varepsilon|^2\right)\md x\ll1\quad(\varepsilon\to0)\,,$$
where $A'=A-\frac{H}{u_\varepsilon^2}\nabla^\bot h_\varepsilon$.\\
From this estimate and the G-L equation satisfied by $\varphi$, we are
able to prove (c.f. \cite[Lemma~6.4]{kach4}) the following estimate
$$\varepsilon\|(\nabla-iA')\varphi\|_{H^1(
  S_1)}\ll1\quad(\varepsilon\to0).$$
Consequently, the trace theorem yields
$$\varepsilon\|n(x)\cdot(\nabla-iA')\varphi\|_{L^2(\partial
  S_1)}\ll1\quad(\varepsilon\to0).$$
Since the functions $h_\varepsilon$ and $u_\varepsilon$ are radial, we
have
$$n(x)\cdot (\nabla-iA')\varphi=n(x)\cdot(\nabla-iA)\varphi\,.$$
Let us also notice  that
$$
\left|n(x)\cdot\left[
\frac{(\nabla-iA)\psi}{\psi}-\frac{\nabla
  u_\varepsilon}{u_\varepsilon}\right]\,\right|
=\frac1{|\varphi|}\left|n(x)\cdot(\nabla-iA')\varphi\right|\,.$$ On
the other hand, since $\frac1{\varepsilon^2}
\int_\Omega(1-|\varphi|^2)^2\,\md x\ll1$, Lemma~\ref{BBH-thm3.3}
yields that $|\varphi|\geq \frac12$ in $\overline\Omega$\,.
Therefore, we deduce that
$$
\varepsilon\left\|n(x)\cdot\left[
\frac{(\nabla-iA)\psi}{\psi}-\frac{\nabla
  u_\varepsilon}{u_\varepsilon}\right]\,\right\|_{L^2(\partial S_1)}
\leq 2\varepsilon
\left\|n(x)\cdot(\nabla-iA')\varphi\right\|_{L^2(\partial S_1)}\ll1\,.$$
Now, invoking Theorem~\ref{mainthm-H=0}, we conclude the result of
Theorem~\ref{thm2}, with $\gamma(a)$ given in (\ref{properties2}).
\hfill$\Box$\\

\subsection{The regime $a>1$}\label{Sec:a>1}
Let us sum up what we know in this case. Let us introduce the
following Ginzburg-Landau functional analyzed in \cite{BBH}
\begin{equation}\label{functional-Serfaty}
H^1(\Omega;\mathbb C)\mapsto
F_\varepsilon(u)=\int_{\Omega}\left(|\nabla
\varphi|^2+\frac1{2\varepsilon^2}(1-|u|^2)^2\right)\,\md x\,.
\end{equation}
Let us also recall the definition of the function
$\xi_\varepsilon:\Omega\rightarrow \mathbb R_-$ introduced in
(\ref{xi-epsilon}). We recall also the set
$\Lambda_\varepsilon=
\{x\in\overline\Omega~:~|\xi_\varepsilon(x)|=\displaystyle
\max_{\overline\Omega}|\xi_\varepsilon|
\,\}$ which governs the location of the vortices of a minimizer of
(\ref{V-EGL}).\\
Now, the result of Lemma~\ref{upperbound-degree} permits to prove
the existence of a constant $\mathcal M>0$ such that (see
\cite[Section~3]{SaSe4})
\begin{equation}\label{Serfaty-Assumption1}
F_\varepsilon(\varphi)< \mathcal M|\ln\varepsilon|\,,\end{equation}
where $\varphi=\frac{\psi}{u_\varepsilon}$ and $(\psi,A)$ always
denote a
minimizer of (\ref{V-EGL}).\\
On the other hand, the result of Theorem~\ref{corollary} states that
\begin{equation}\label{Serfaty-Assumption2}
\Lambda_\varepsilon=\{0\}\,,\quad \xi_\varepsilon''(0)>0\,.
\end{equation}
The estimate (\ref{Serfaty-Assumption1}) is  the basis on which the
analysis in \cite{Se1} is build-up. It permits to prove an
expression of the first critical field:
$$H_{C_1}=k_\varepsilon|\ln\varepsilon|+k_{1,\varepsilon}\,,$$
where $k_\varepsilon$ is given by (\ref{k-epsilon}) and
$k_{1,\varepsilon}=\mathcal O(1)$. If $H_{C_1}+k<H<H_{C_1}+\mathcal
O(1)$,
$k>0$,  then  a
minimizer $(\psi,A)$ of (\ref{V-EGL}) has a finite number of
vortices, each of degree $1$, and localized near the center of the
disc $\Omega=D(0,1)$. Furthermore, it is proved that if more than
one vortex exists, distinct vortices will tend, after normalization,
to distinct points in
$\mathbb R^2$.\\
The results in (\ref{Serfaty-Assumption2}) are the basis to build-up
the analysis of \cite{Se2}, which permits to obtain a sequence of
critical fields. We point out that in order to adapt the analysis of
\cite{Se2}, we need to remember that in every compact subset $K$ of
$D(0,R)$,  the function $u_\varepsilon$ converges to $1$
exponentially fast in $C^2(K)$.\\
We define the normalized energy~:
\begin{equation}\label{normalizedenergy}
w_{\varepsilon,n}: \mathbb R^n\ni(x_1,x_2,\cdots, x_n)\mapsto
-2\pi\sum_{i\not=j}\ln|x_i-x_j|+2\pi
\xi_\varepsilon''(0)\sum_{i=1}^n|x_i|^2\,.
\end{equation}
The analysis of \cite{Se2} yields that, if the magnetic field
satisfies
$$H=k_\varepsilon\left(|\ln\varepsilon|+\delta\ln|\ln\varepsilon|\right),\quad
n-1<\delta<n,\quad n\in\mathbb N\,,$$ then a minimizer $(\psi,A)$ of
(\ref{V-EGL}) has $n$-vortices $(x_i(\varepsilon))_{i=1}^n$, each of
degree $1$, and such that, when putting $\widetilde
x_i(\varepsilon)=x_i(\varepsilon)\sqrt{H}\,$, then  the
configuration $(\widetilde x_i(\varepsilon))_{i=1}^n$ is localized
near a minimizer of the renormalized  energy $w_{\varepsilon,n}$.
Furthermore, the following expansion of the energy holds as
$\varepsilon\to0$\,:
\begin{eqnarray*}
\mathcal
F_{\varepsilon,H}(\varphi,A)&=&H^2J_0(\varepsilon)
-2\pi\,n\left(|\ln\varepsilon|-\frac{H}{k_\varepsilon}\right)
+\pi(n^2-n)\ln H\\
&&+w_{\varepsilon,n}\left(\widetilde
x_1(\varepsilon),\cdots,\widetilde
x_n(\varepsilon)\right)+Q_n+o(1)\,,
\end{eqnarray*}
where $Q_n$ is an explicit constant
depending only on $n$.

\end{document}